\documentclass[12pt]{article}
\usepackage{amsmath,amssymb,amsthm,esint,dsfont,tikz}
\usepackage[english]{babel}
\usepackage[margin=1in]{geometry}
\usepackage{csquotes}
\usepackage{setspace}
\usepackage[colorlinks, pdfpagelabels, pdfstartview = FitH, bookmarksopen=true, bookmarksnumbered = true, linkcolor = black,plainpages = false, hypertexnames = false, citecolor = black]{hyperref}

\usepackage{color}

\newcommand{\answ}[1]{#1}
\newcommand{\blue}[1]{#1}

\newcommand{\dth}[1]{#1}

\renewcommand{\H}{\mathcal{H}}
\newcommand{\exe}{E_{\xi\dth{,\eta}}}

\newcommand{\beqn}{\begin{equation}}
\newcommand{\eeqn}{\end{equation}}

\newcommand{\Om}{\Omega}
\newcommand{\sym}{\mathcal{S}_0}

\newcommand{\grad}{\nabla}
\newcommand{\tilQ}{\widetilde Q}
\newcommand{\Qm}{Q_{\xi\dth{,\eta}}}
\newcommand{\C}{\mathcal{C}}

\renewcommand{\sp}{{\mathbb{S}^2}}
\newcommand{\Qc}{Q_{\xi\dth{,\eta}}^{h,\varepsilon}}
\newcommand{\Qcon}{\overline{Q_{\xi\dth{,\eta}}}}
\renewcommand{\S}{\mathbb{S}}

\newcommand{\abs}[1]{{\left\vert #1\right\vert}}
\newcommand{\e}{\varepsilon}
\newcommand{\R}{{\mathbb R}}

\DeclareMathOperator{\supp}{supp}
\DeclareMathOperator{\tr}{tr}
\renewcommand{\L}{\mathcal{L}}
\newcommand{\Qvb}{\overline{Q_v}}
\newcommand{\Qib}{\overline{Q_\infty}}
\newcommand{\Qmb}{\overline{Q_{\mathrm{mix}}}}
\newcommand{\Qbb}{\overline{Q_v}}

\renewcommand{\hat}{\widehat}

\newtheorem{theorem}{Theorem}[section]
\newtheorem{proposition}[theorem]{Proposition}
\newtheorem{lemma}[theorem]{Lemma}

\theoremstyle{definition}
\newtheorem{remark}[theorem]{Remark}

\usepackage[capitalise]{cleveref}
\usepackage{authblk}

\title{Spherical Particle in Nematic Liquid Crystal with a Magnetic Field and Planar Anchoring}
\date{\today}

\begin{document}

\author[1]{Lia Bronsard\footnote{bronsard@mcmaster.ca }}
\author[1]{Dean Louizos\footnote{louizod@mcmaster.ca}}
\author[1]{Dominik Stantejsky\footnote{stantejd@mcmaster.ca}}
\affil[1]{Department of Mathematics and Statistics, McMaster University, Hamilton, ON L8S 4L8 Canada}

\maketitle

\begin{abstract}
We study minimizers of the Landau-de Gennes energy in $\mathbb{R}^3\setminus B_1(0)$ with external magnetic field in the large particle limit.
We impose strong tangential anchoring and uniaxiality of the $Q-$tensor on the boundary.
We derive a lower bound for the energy in terms of the boundary condition and show in the extreme cases \answ{of strong and weak magnetic field strength} that the longitudinal director field is energy minimizing, indicating the presence of two half-point defects, so called \emph{boojums}, at two opposite points of the sphere.
Using a recovery sequence, we show that the energy bound is optimal \answ{in these extreme cases.} \\
\linebreak
\textbf{Keywords:} Singular limit, tangential anchoring, boojums, Landau-de Gennes energy.
\linebreak
\textbf{MSC2020:} 
49J45, 
35J50, 
49S05, 
76A15. 
\end{abstract}

\section{Introduction}
\label{sec:intro}

In this paper we study the effect of an external magnetic field on the type of defects that can be observed near
a spherical colloidal particle immersed in nematic liquid crystal with tangential anchoring. This problem has been studied in the case of normal Dirichlet boundary conditions in \cite{weakanchor,colloidphys,abl,ACS2021,fukudaetal04,fukudayokoyama06} and in this article we are interested in the effect of tangential boundary conditions. 
Such tangential anchoring is experimentally observed and investigated \cite{Liu2013,Senyuk2021,Shi2005,Tas2012} but very little mathematical analysis has been done for these boundary conditions and new analytical questions arise in this context \cite{abc,ABv,bcs}. 

An important feature in experiments is the occurrence of defects leading to a challenging but interesting  analysis of possible singular structures \answ{such as point or line defects, on the boundary or in the interior of the liquid crystal.}
Various mathematical models have been proposed to describe nematic liquid crystals.
In this article, we use the Landau-de-Gennes model in which the orientation of the elongated molecules is represented by a $Q$-tensor, namely by a traceless symmetric $3\times 3$ matrix, accounting for the alignment of the molecules through their eigenvectors and eigenvalues. 
Another continuum model is the Oseen-Frank model, in which the orientation is represented by a unit director $n\in\mathbb S^2$, and the $Q$-tensors can be thought of as a relaxation of this uniaxial constraint. 
The Oseen-Frank model can be recovered from the Landau-de Gennes energy in the limit as the \answ{nematic coherence length, i.e.\ the typical size of a defect core, becomes small relative to the size of the colloidal particle, see \cite{Gartland} for more details on this limit}.
This convergence has recently produced new powerful mathematical analysis \cite{baumanparkphillips12,
canevari2d,
canevari3d,
singperturb,DiRoSlZa,
golovatymontero14,
majumdarzarnescu10,NguZa}.  
Compared to the Oseen-Frank model, the additional degrees of freedom given by $Q$-tensors allow for a much finer description of the defect cores \answ{e.g.\ when the liquid crystal molecules have more than one preferred director. This
translates to $Q$ having three distinct eigenvalues (``biaxiality'') rather than two equal eigenvalues (``uniaxiality''), and the $Q$-tensor approach allows to bypass problems e.g.\ due to orientability \cite{BaZar}.} 
The nonlinear analysis of defect cores has recently attracted much attention \cite{canevari2d,CaRaMa,biaxialescape,difrattaetal16,INSZuniqhedgehog,
INSZstabhedgehog,INSZstab2d,INSZinstab2d}.
Further, when colloidal particles are immersed into a nematic liquid crystal, they induce distortions in the nematic alignment which may generate additional defects, leading to many potential applications related e.g.\ to the detection of these inclusions or to structure formation created by defect interactions \cite{stark01,Shi2005}. 
The mathematical analysis of such phenomena is very challenging and we are interested in the most fundamental situation: a single spherical particle surrounded by a liquid crystal. 

As opposed to the previous results mentioned above, we will study planar anchoring at the particle surface.\answ{ This corresponds to the case when}
the liquid crystal molecules \answ{are parallel to the two dimensional surface, hence} leaving an additional degree of freedom in the tangent plane at the surface. 
This degree of freedom leads to interesting new mathematically challenging questions, including the choice of preferred direction in the tangent space. 
We are able to show in certain regimes that the longitudinal direction is the preferred one, and leave open as conjecture that it is always the case.

In the case of the Dirichlet boundary conditions known as radial anchoring, when the liquid crystal molecules align perpendicularly to the surface of the particle, it has been shown \cite{colloid} that in the small particle limit, a \emph{Saturn ring} defect occurs. 
Further, when an external field is imposed,  it was shown that a \emph{Saturn ring} defect can still occur for large particles \cite{abl, ACS2021}. 
In this article, we are interested in studying the effect of a uniform magnetic field in the sample on the type of defect expected when tangential anchoring is imposed for large particles. 
We will be using the simplified model from \cite{fukudaetal04,fukudayokoyama06}  where the presence of the external field is modelled by adding a symmetry-breaking term to the energy (favoring alignment along the field), multiplied by a parameter accounting for the intensity of the field.

After adequate non-dimensionalization \cite{fukudaetal04} we are left with two parameters $\xi,\eta>0$ which represent, in units of the particle radius, the coherence lengths for nematic alignment and magnetic coherence length. 
In these units, the colloid particle is represented by the closed ball of radius one $B:=B_1(0)\subset\R^3$, so that the liquid crystal is contained in the domain  $\Omega=\R^3\setminus B$. 
The Landau-de Gennes energy used in \cite{abl,ACS2021, fukudaetal04,fukudayokoyama06} is given  by
\begin{align}\label{energy}
E_\answ{\xi,\eta}(Q)
\ &= \
\int_{\Omega} \left( \frac 12\abs{\nabla Q}^2 +\frac{1}{\xi^2}\left[f(Q) + h^2 g(Q) \right]\right) \, dx  \notag\\
\ &= \
\int_{\Omega} \left( \frac 12\abs{\nabla Q}^2 +\frac{1}{\xi^2}f(Q) + \frac{1}{\eta^2} g(Q) \right) \, dx 
\, ,
\end{align}
where $\eta=\frac\xi h$ and $h$ is the imposed magnetic field strength.
The map $Q$ takes values into the space $\sym$ of $3\times 3$ symmetric matrices with zero traces.
The nematic potential at a given constant temperature can be assumed to be given, for non negative constants $a,b,c,C$, by
\begin{equation}\label{eq:def:f}
f(Q)
\ = \
C - \frac a2\abs{Q}^2 - \frac b3\tr(Q^3) + \frac c4 \abs{Q}^4
\, .
\end{equation}
\answ{For specified values of $a,b,c\geq 0$, one can choose $C\geq0$ such that} $f$ satisfies
\begin{equation*}
f(Q) \ \geq \ 0
\text{ with equality iff }
Q \ = \ \answ{s_*\Big(}n\otimes n -\frac 13 I\answ{\Big)}
\text{ for some }
n\in\mathbb S^2
\, ,
\end{equation*}
\answ{where $s_*:=\frac{1}{4c}(b+\sqrt{b^2+24ac})$. 
For simplicity, we will choose $a=1$ and $b=3=c$ so that $s_*=1$, however, our analysis does not require this assumption.}
The symmetry-breaking potential $g(Q)$ is given by
\begin{equation*}
g(Q)
\ = \
\sqrt{\frac 23} - \frac{Q_{33}}{\abs{Q}}
\, .
\end{equation*}
It breaks symmetry in the sense that the rotations $R\in SO(3)$ which satisfy $g(R^\top Q R)=g(Q)$ for all $Q\in\sym$ must have $e_3$ as an eigenvector, while $f(R^\top Q R)=f(Q)$ for all $R\in SO(3)$ and $Q\in \sym$. 
Its specific form is chosen so that
\begin{equation*}
g(Q)
\ = \
\sqrt{\frac32}(1-n_3^2)
\quad\text{for }Q \ = \ n\otimes n-\frac 13 I
\, ,
\end{equation*}
and $g(Q)$ is invariant under multiplication of $Q$ by a positive constant \cite{fukudaetal04}.
This potential satisfies
\begin{equation*}
g(Q)
\ \geq \ 0
\text{ with equality iff }
Q \ = \ 
\mu\left(\mathbf e_3\otimes \mathbf e_3 -\frac 13 \right)I
\text{ for some }
\mu > 0
\, .
\end{equation*}
The form of the potential $g(Q)$ that models the influence of the external magnetic field is the result of simplifying assumptions as detailed in \cite{fukudaetal04,fukudayokoyama06}.
In particular, the field is taken to be constant throughout the liquid crystal sample, which is more realistic to assume in the case of a magnetic, compared to an electric field. 
We do not question the physical validity of these assumptions, but focus on this model where the external field introduces a symmetry-breaking effect at some additional length scale.

For $h>0$ the full potential $f(Q)+h^2 g(Q)$ is minimized exactly at $Q=Q_\infty$, where
\begin{equation*}
Q_\infty
\ = \
\mathbf e_3 \otimes \mathbf e_3 -\frac 13 I
\, .
\end{equation*}
Moreover, it is easily checked that
\begin{equation*}
f(Q) + h^2 g(Q)
\ \geq \
C(h)\abs{Q-Q_\infty}^2
\, ,
\end{equation*}
for some constant $C(h)>0$. 
This ensures that the energy is coercive on the affine space $Q_\infty + H^1(\Omega;\sym)$.
The anchoring at the particle surface is assumed to be tangential:
\begin{equation*}
Q
\ = \ 
Q_b
\ := \ 
v\otimes v -\frac 13 I
\quad\text{on }\partial \Omega,
\quad  v\in \mathbb{S}^2, \quad \text{ with }
\quad  v\answ{(\omega)}\cdot \nu\answ{(\omega)}=0\answ{\quad\text{for all $\omega\in\partial \Om$}}
\, ,
\end{equation*}\answ{where $v:\partial \Om\to\S^2$ is measurable and} $\nu$ is the outward unit normal to $\answ{\partial}\Om$.
A particularly important choice is the tangent vector $v=e_\phi$, corresponding to the longitudinal direction on the sphere and we denote $Q_\phi := e_\phi\otimes e_\phi -\frac 13 I$.
Letting $\mathcal H$ be the space
\begin{equation*}
\mathcal H 
\ = \
\left\lbrace Q\in Q_\infty + H^1(\Omega;\mathcal S_0)\colon Q=Q_b\text{ on }\partial\Omega\right\rbrace
\, ,
\end{equation*}
the coercivity of the energy ensures existence of a minimizer in $\mathcal H$ for any $\xi,\eta>0$.

The combined potentials
\begin{equation*}
\frac{1}{\xi^2}f(Q)+\frac{1}{\eta^2}g(Q) = \frac{1}{\eta^2}\left(\frac{\eta^2}{\xi^2}f(Q)+g(Q)\right),
\end{equation*}
is minimized at $Q=Q_\infty$. 
As $\eta\to 0$ this forces a minimizing configuration to be very close to $Q_\infty$. 
Note that the boundary data $Q_b$ satisfies $f(Q_b)\equiv 0$ but not necessarily $g(Q_b) = 0$. 
Not surprisingly, as in \cite{abl,ACS2021}, deformations concentrate in a boundary layer of size $\eta$, where a one-dimensional transition takes place according to the energy
\begin{equation}\label{Flambda}
F_\lambda(Q)=\int_1^\infty \left[\frac 12  \abs{\frac{d Q}{d r}}^2 +\lambda^2 f(Q)+g(Q)\right] \, dr,
\end{equation}
\answ{where $\lambda\in [0,\infty]$. 
Later on, $\lambda$ is identified with the limit of $\frac{\eta}{\xi}$ as $\xi,\eta\to 0$. 
Therefore, $\lambda$ can be seen as a relative strength of the magnetic field w.r.t\ the elastic forces, a large (resp.\ small) $\lambda$ corresponding to a weak (resp.\ strong) magnetic field.
The energy \eqref{Flambda} is well}
defined for $Q\in Q_\infty + H^1((1,\infty);\mathcal S_0)$.
For $\lambda=\infty$ the formula in \eqref{Flambda} should be understood as
\begin{equation*}
F_\infty(Q)=
\left\lbrace
\begin{aligned}
&\int_1^\infty \left[\frac 12  \abs{\frac{d Q}{d r}}^2 + g(Q)\right] \, dr &&\text{ if }f(Q)=0\text{ a.e.},\\
&+\infty && \text{ otherwise.}
\end{aligned}
\right.
\end{equation*}
In other words, $F_\infty$ is finite for maps $Q\in Q_\infty + H^1((1,\infty);\mathcal S_0)$ which satisfy $Q=n\otimes n -\frac13 I$ for some measurable map $n\colon (1,\infty)\to\mathbb S^2$.

Obviously, the cases $\lambda\in [0,\infty)$ and $\lambda=\infty$ are quite different and they require separate treatments, but in both cases we obtain for the energy $\exe$ 
of a minimizer $\Qm$ the asymptotics 
\begin{equation*}
\eta\: \exe(\Qm;\Omega)
\ = \
\int_{\mathbb S^2}D_\lambda(Q_b(\omega))\,d\mathcal H^2(\omega) 
+ o(1)
\qquad\text{as }
\xi\dth{,\eta}\to0
\, ,
\end{equation*}
where \dth{$\exe$ is from \eqref{energy} and}
\begin{equation}\label{eq:Dlambda}
D_\lambda(Q_0)
\ = \
\min\left\lbrace F_\lambda(Q)\colon Q\in Q_\infty +H^1((1,\infty);\mathcal S_0),\,Q(1)=Q_0\right\rbrace
\, .
\end{equation}
The existence of a minimizer of $F_\lambda$ which attains $D_\lambda(Q_b(\omega))$ for any $\omega\in\mathbb S^2$ follows from the direct method.  
In section~\ref{ss:upII} we will exploit the observation of Sternberg \cite{St91} that the heteroclinic connections which minimize $F_\lambda$ represent geodesics for a degenerate metric.
In the case $\lambda=\infty$ this enables us to obtain an exact value for the limiting energy in terms of the angles $\theta\in[0,2\pi)$ and $\phi\in [0,\pi]$ \answ{which parameterize the unit sphere}:
\begin{equation*}
D_\infty(Q_b(\theta,\phi))
\ = \
\kappa(1-|\sin\theta|), 
\quad \text{and}\quad 
\lim_{\xi\dth{,\eta}\to0} \eta\, \exe(\Qm;\Om) 
\ = \ \pi(4-\pi)\kappa
\, ,
\end{equation*}
where $\kappa:=\sqrt[4]{24}$.
Note that in \cite{abl} the corresponding minimal energy for radial anchoring has been computed to be $D_\infty(Q_b(\theta,\phi)) = \kappa(1-|\cos\theta|)$ and the asymptotic minimal value is $2\pi\kappa$.

More specifically, we obtain local asymptotics in conical subdomains of $\Omega$. 
For $U\subset\mathbb S^2$ we denote by $\mathcal C(U)$ the cone
\begin{equation*}
\mathcal C(U)
\ := \ 
\left\lbrace t\omega\colon t>1,\,\omega\in U\right\rbrace
\, ,
\end{equation*}
and prove
\begin{theorem}[Lower bound]\label{thm:lower}
    Let $\Qm$ minimize $\exe$ with $\Qm=Q_b$ on $\sp$. If
    \begin{equation*}
        \frac{\eta}{\xi}\to\lambda\in [0,\infty],\quad\text{as}\quad\xi\dth{,\eta}\to0,
    \end{equation*}
    then for any measureable set $U\subset \sp$,
    \begin{equation*}
        \liminf_{\xi\dth{,\eta}\to0}\: \eta\:  \exe(\Qm\dth{;\, }\C(U))
        \ \geq \
        \int_U D_\lambda(Q_b(\omega))\ \answ{d\H^2(\omega)}.
    \end{equation*}
\end{theorem}

The lower bound in Theorem~\ref{thm:lower} follows from an elementary rescaling and the properties of $\lambda\mapsto D_\lambda$ as in \cite{abl} with only minor modifications.
Note that this lower bound is valid for any admissible $Q_b$.
So a natural question is, for which choice of $Q_b$ the energy density $D_\lambda(Q_b(\omega))$ attains its minimum. 
\answ{We conjecture that for all $\lambda\in [0,\infty]$, the function $Q_\phi$ minimizes $D_\lambda(Q_b(\omega))$ (see Section 3).}
We are able to answer this question in the following cases:

\begin{theorem}[Minimizing boundary conditions]\label{thm:optim}
    If
    \begin{equation*}
        \frac{\eta}{\xi}\to\lambda\in\{0,\infty\},\quad\text{as}\quad\xi\dth{,\eta}\to0,
    \end{equation*}
    then for any point $\omega\in \sp$,
    \begin{equation*}
        D_\lambda(Q_b(\omega))
        \ \geq \
        D_\lambda(Q_\phi(\omega))
        \, ,
    \end{equation*}
    for any admissible $Q_b$.
\end{theorem}

For $\lambda=\infty$, the minimizers $Q$ of $F_\lambda$ are explicitly known, see Lemma \ref{def:n_min}.
If $\lambda<\infty$, the precise profile is harder to obtain.
Nevertheless, in Lemma \ref{lem:lambda:zero:shape-minimizers} we obtain a characterization of the minimizer for $\lambda=0$ in terms of a solution to a scalar ODE.
In the intermediary case $\lambda\in (0,\infty)$, we are still able to simplify the problem to a system of three coupled ODEs, see Proposition \ref{prop:lambda:finite}.

\answ{Assuming the boundary condition $Q_b=Q_\phi$, we can }construct an upper bound that matches the lower bound in Theorem~\ref{thm:lower} \answ{for $Q_b=Q_\phi$}:

\begin{theorem}[Recovery sequence]\label{thm:upper}
    If
    \begin{equation*}
        \frac{\eta}{\xi}\to\lambda\in [0,\infty],\quad\text{as}\quad\xi\dth{,\eta}\to0,
    \end{equation*}
    then there exists $\Qm$ with $\Qm=Q_\answ{\phi}$ on $\sp$ such that for any measureable set $U\subset \sp$,
    \begin{equation*}
        \limsup_{\xi\dth{,\eta}\to0}\: \eta\:  \exe(\Qm\dth{;}\, \C(U))
        \ \leq \
        \int_U D_\lambda(Q_\answ{\phi}(\omega))\ \answ{d\H^2(\omega)}.
    \end{equation*}
\end{theorem}

\answ{We remark that in the case $\lambda\in\{0,\infty\}$ the asymptotic energy of minimizers is completely determined by Theorem~\ref{thm:lower} and \ref{thm:upper} since, according to Theorem~\ref{thm:optim}, $Q_\phi$ is the minimizing boundary condition.}

\answ{Although the} result \answ{in Theorem \ref{thm:upper}} is analogue to \cite{abl}, the construction necessary to obtain Theorem~\ref{thm:upper} is very different.
The main distinction is that the singular structure appearing for $Q_b=Q_\phi$ is no longer a line defect (with isotropic or biaxial core), but two half-point defects (so called \emph{boojums}).
Those defects are located at opposite poles and can be constructed to be uniaxial everywhere.
To realize this recovery sequence e.g.\ for $\lambda=\infty$, it is necessary that all construction are done within the uniaxial manifold.
In particular, the interpolations between point defects and optimal profile achieving the minimal energy density $D_\infty(Q_b)$ 
\answ{cannot be} a direct linear interpolation of the $Q-$tensor components.
\answ{Instead, we interpolate the angles between the director fields of the boundary condition, of the defect region, of the optimal profile and $e_3$, and then use this angle to define a uniaxial $Q-$tensor, see the proof of Proposition~\ref{prop:up:inf}.
} 
\answ{
\begin{remark}\label{rmk}
    The method to construct the recovery sequence presented in the proof of Theorem \ref{thm:upper} can be directly adapted to the case of an equivariant Lipschitz boundary condition $Q_b$ which agrees locally with $Q_\phi$ near the poles, and this is made precise in the proof by using $Q_\phi$ only when necessary. 
    The construction introduced in \cite{BLS2} could be used to further 
    enlarge the class of boundary conditions for which an upper bound (matching the lower bound from Theorem~\ref{thm:lower}) holds,
    namely only requiring  that $Q_b$ is Lipschitz away from finitely many isolated points.
\end{remark}
} 

\paragraph{Organization of the paper.} 
In Section \ref{sec:lower} we prove the lower bound. 
We then discuss the optimal choice of $Q_b$ for the limit energy in Section \ref{sec:optim-bc}.
In Section \ref{sec:upper} we concentrate on the upper bound, considering the case $\lambda\in [0,\infty)$ in Subsection \ref{ss:upI} and $\lambda=\infty$ in Subsection \ref{ss:upII}. 

\paragraph{Acknowledgements:} 
LB is supported by an NSERC (Canada) Discovery Grants.

\section{Lower Bound}
\label{sec:lower}

We begin by considering a general planar boundary condition $Q_b=v\otimes v-\frac{1}{3}I$, where $v:\S^2\to\S^2$ is such that for a.e.\ $\omega\in\sp$, $\answ{\nu}(\omega)\cdot v(\omega)=0$. We also impose that $v$ is smooth a.e.\ on $\S^2$ and is discontinuous only at finitely many isolated points. 
We recall that for a measurable set $U\subset \sp$ we defined the cone $\C(U):=\{r\omega:\omega\in\answ{U},\ 1<r<\infty\}$. 

In this section we prove Theorem \ref{thm:lower}.
The proof is very similar to \cite{abl} where the case $Q_b=e_r\otimes e_r-\frac{1}{3}I$ was considered.
For the convenience of the reader we recall the main steps here:

We first change to spherical coordinates by $x=r\omega$, $(r,\omega)\in (1,\infty)\times\mathbb S^2$, then define
\begin{equation*}
    \widetilde Q(\tilde r,\omega):=\Qm(1+\eta(\tilde r -1),\omega),
\end{equation*}
where $r=1+\eta(\tilde{r}-1)$. Applying this change of variables, we have
\begin{align}
    \eta \exe(\Qm,\mathcal C(U))&
    =\eta\int_U \int_1^\infty  \left[\frac{1}{2} \abs{\grad \Qm}^2+\frac{1}{\xi^2}f(\Qm)+\frac{1}{\eta^2}g(\Qm)
    \right]\,r^2\, dr \, \answ{d\H^2(\omega)}\nonumber \\
    &\geq \int_U \int_1^\infty  \left[\frac{1}{2} \abs{\frac{\partial\tilQ}{\partial\tilde r}}^2+\frac{\eta^2}{\xi^2}f(\tilQ)+g(\tilQ)
    \right]\, d\tilde r \, \answ{d\H^2(\omega)} \label{prop:lower:eq}
\end{align}
Then, for $\lambda\in (0,\infty]$, using the definitions \eqref{Flambda} and \eqref{eq:Dlambda}, we see that
\begin{equation*}
    \eta \exe(\Qm\dth{;\, }\C(U))\geq \int_U F_\frac{\eta}{\xi}(\tilQ\dth{(\cdot,\omega)})\, \answ{d\H^2(\omega)}\geq \int_U D_{\frac\eta\xi}(Q_b(\omega))\,\answ{d\H^2(\omega)},
\end{equation*}
We now use Lemma 2.2 from \cite{abl} which states that for any $Q_0\in\sym$ and $\lambda\in(0,\infty]$,
\[D_\lambda(Q_0)=\lim_{\mu\to\lambda}D_\mu(Q_0).\]
This result along with Fatou's lemma gives the desired lower bound for for $\lambda\in (0,\infty]$
\begin{equation*}
    \liminf_{\xi\dth{,\eta}\to0}\eta \exe(\Qm\dth{;\, }\C(U))\geq \int_U D_\lambda(Q_b(\omega))\, \answ{d\H^2(\omega)}.
\end{equation*}

If $\lambda=0$, then one can directly estimate \eqref{prop:lower:eq} by
\begin{align*}
    \liminf_{\xi\dth{,\eta}\to0}\eta \exe(\Qm\dth{;\, }\C(U))
    \ \geq \ 
    \int_U F_0(\widetilde{Q}\dth{(\cdot,\omega)})\, \answ{d\H^2(\omega)}
    \ \geq \
    \int_U D_0(Q_b(\omega))\, \answ{d\H^2(\omega)}
    \, .
\end{align*}

We note that this bound depends explicitly on the choice of function $v$ in the definition of $Q_b$, so we can find a smaller lower bound of the same form by choosing a specific $v$ and this will be discussed in the next section.

\section{Optimal Boundary Conditions}
\label{sec:optim-bc}

The goal of this section is to prove that $D_\lambda(Q_b)$ is minimized when we choose $Q_b=e_\phi\otimes e_\phi-\frac{1}{3}I$. 
Intuitively this is expected as $(\mp) e_\phi$ is the unit tangent vector closest to $(\pm) e_3$ at each point on $\sp\setminus\{\pm e_3\}$. 
The three cases of $\lambda=\infty$, $\lambda=0$ and $\lambda\in(0,\infty)$ are vastly different, so we consider them separately, beginning with $\lambda=\infty$.

\begin{proposition}[The case $\lambda=\infty$]\label{prop:min:inf}
    Let $\omega\in \sp$ be fixed, define $Q_\phi=e_\phi\otimes e_\phi-\frac{1}{3}I$ and $Q_v=v\otimes v-\frac{1}{3}I$, where $v\in\S^2$ with $v\cdot \omega=0$, then
    \[D_\infty(Q_\phi)\leq D_\infty(Q_v)\]
    with equality if and only if $v=\pm e_\phi$.
\end{proposition}
Due to the symmetry of the problem, it will suffice to prove this for $\omega\in \sp$ with $0\leq\omega_3<1$, although the case where $\omega_3=0$ is trivial. It will also make sense to only consider vectors $v\in\sp$ such that $v_1<0$, $v_2=0$ and $0\leq v_3<1$. 
This will be explained further in the proof of Proposition \ref{prop:min:inf}, but for now we consider the following related lemma which deals with the quantity $G_\infty$ defined to be,
\begin{align*}
    G_\infty(n)
    \ := \
    \int_0^\infty\abs{\frac{\partial n}{\partial t}}^2+g(n)\ dt
    \, ,
\end{align*}
where $g(n):=\sqrt{\frac32}(1-n_3^2)=g(Q)$ for $Q=n\otimes n - \frac13 I$.

This lemma is a reformulation of \cite[Lemma 3.4]{abl} and provides an explicit description of the minimisers of $G_\infty$ among functions $n\in H^1((0,\infty),\mathbb{R}^3)+e_3$ with $|n|=1$ and a given initial condition.

\begin{lemma}\label{def:n_min}
    Let $v\in\sp$ with $v_1<0$, $v_2=0$ and $0\leq v_3<1$, then if $n$ minimizes $G_\infty$ with $n(0)=v$, $n$ is of the form, $n=(-\sqrt{1-n_3^2},0,n_3)$, where
    \begin{equation}
        n_3(t,\varphi)=\frac{A(\varphi)-e^{-\kappa t}}{A(\varphi)+e^{-\kappa t}},\quad\text{and}\quad A(\varphi)=\frac{1+\cos\varphi}{1-\cos\varphi},\label{n_min:eq}
    \end{equation}
    with $\varphi$ being the angle between $v$ and $e_3$ and with $\kappa=\sqrt[4]{24}$. Further, there exists a constant $C>0$ independent of $\varphi$ such that
    \begin{equation}\label{n_min:bound}
        \abs{\frac{\partial n}{\partial t}}^2,\abs{\frac{\partial n}{\partial \varphi}}^2,|n_1|^2\leq Ce^{-\kappa t}.
    \end{equation}
\end{lemma}

Using this Lemma, we can proceed with the proof of Proposition \ref{prop:min:inf}.

\begin{proof}[Proof of Proposition \ref{prop:min:inf}]
We begin by treating the case $\omega=\pm e_3$, in which the rotational symmetry of the problem implies that all unit tangent vectors $v$ have equal energy $D_\infty(Q_v)$.

Next, we consider $\omega\in \sp$ with $\omega_3=0$, then $Q_\phi=Q_\infty$, so 
\[D_\infty(Q_\phi)=0\leq D_\infty(Q_v).\]
We use the symmetry of the setup to reduce the problem to the case where $0<\omega_3<1$, $v_1<0$, $v_2=0$ and $0\leq v_3<1$.

First, since $Q_v=Q_{-v}$ for any $v\in\sp$, it holds $D_\infty(Q_v)=D_\infty(Q_{-v})$\dth{, so we can consider $0\leq v_3<1$}.
If $D_\infty(Q_\phi)\leq D_\infty(Q_v)$ at $\omega$, the same holds at $-\omega$, so we can impose without loss of generality that $0<\omega_3<1$. 
Now let $v=(v_1,v_2,v_3)$.
Then there exists $\theta_0$ and a rotation about the $x_3$-axis $R_{\theta_0}$ such that for $u=R_{\theta_0}v$, we have $u_1<0$ and $u_2=0$.
If $n$ minimizes $G_\infty$ subject to $n(0)=u$, then $\tilde n=R_{-\theta_0}n$ minimizes $G_\infty$ subject to $\tilde n(0)=v$ and \[D_\infty(Q_v)=G_\infty(\tilde n)=G_\infty(n).\]
In this setup, $D_\infty(Q_v)=G_\infty(n)$ where $n$ is as defined in Lemma \ref{def:n_min}, so we can compute $D_\infty(Q_v)$ explicitly.
\begin{align*}
    D_\infty(Q_v)=\int_0^\infty\abs{\frac{\partial n}{\partial t}}^2+g(n)\ dt=\int_0^\infty 2\abs{\frac{\partial n}{\partial t}}\sqrt{g(n)}\ dt
    \, ,
\end{align*}
where $n$ was chosen specifically to satisfy this equality. 
By a direct computation, we get that
\begin{align*}
    D_\infty(Q_v) 
    \ = \ 
    \kappa\int_0^\infty\abs{\frac{\partial n_3}{\partial t}}\ dt=\kappa(1-v_3)
    \, .
\end{align*}
We note that if $v\cdot\omega=0$, then $v_3$ is maximized for $v=- e_\phi$, and thus $D_\infty(Q_v)$ is minimized by $Q_\phi$.
\end{proof}

When dealing with $\lambda\in[0,\infty)$ we will decompose $Q$ into two components, $N=|Q|$ and $\hat{Q}=Q/|Q|$. Subsequently it will make sense to split the energy $F_\lambda$ (recall the definition  \eqref{Flambda}) into two components $H(N)$ and $K_N(\hat{Q})$ defined as,
\begin{align}
    F_\lambda(Q)&=\int_1^\infty\frac{1}{2}|\partial_r N|^2+\frac{1}{2}N^2|\partial_r\hat{Q}|^2+\lambda^2\Big(\frac{2}{9}-\frac{1}{2}N^2-N^3\tr(\hat{Q}^3)+\frac{3}{4}N^4\Big)+g(\hat{Q})\, dr\notag\\
    &=\int_1^\infty\frac{1}{2}|\partial_r N|^2+\lambda^2\Big(-\frac{1}{2}N^2+\frac{3}{4}N^4\Big)\, dr\ \notag\\
    &\quad\qquad+\ \int_1^\infty \frac{1}{2}N^2|\partial_r\hat{Q}|^2+\lambda^2\Big(\frac{2}{9}-N^3\tr(\hat{Q}^3)\Big)+g(\hat{Q})\, dr\notag\\
    &=:H(N)+K_N(\hat{Q})
    \, .\label{H-Kn}
\end{align}
We start with the simpler case $\lambda=0$.

\begin{proposition}[The case $\lambda=0$]\label{prop:lambda:zero}
    Let $\omega\in\sp$ be fixed, then
    \begin{equation*}
        D_0(Q_\phi)\leq D_0(Q_v)
    \end{equation*}
    for $Q_\phi$ and $Q_v$ defined in Proposition \ref{prop:min:inf}.
\end{proposition}

In order to prove this proposition, we first need the following lemma.

\begin{figure}[hbt!]
\begin{center}
\includegraphics[scale=1.5]{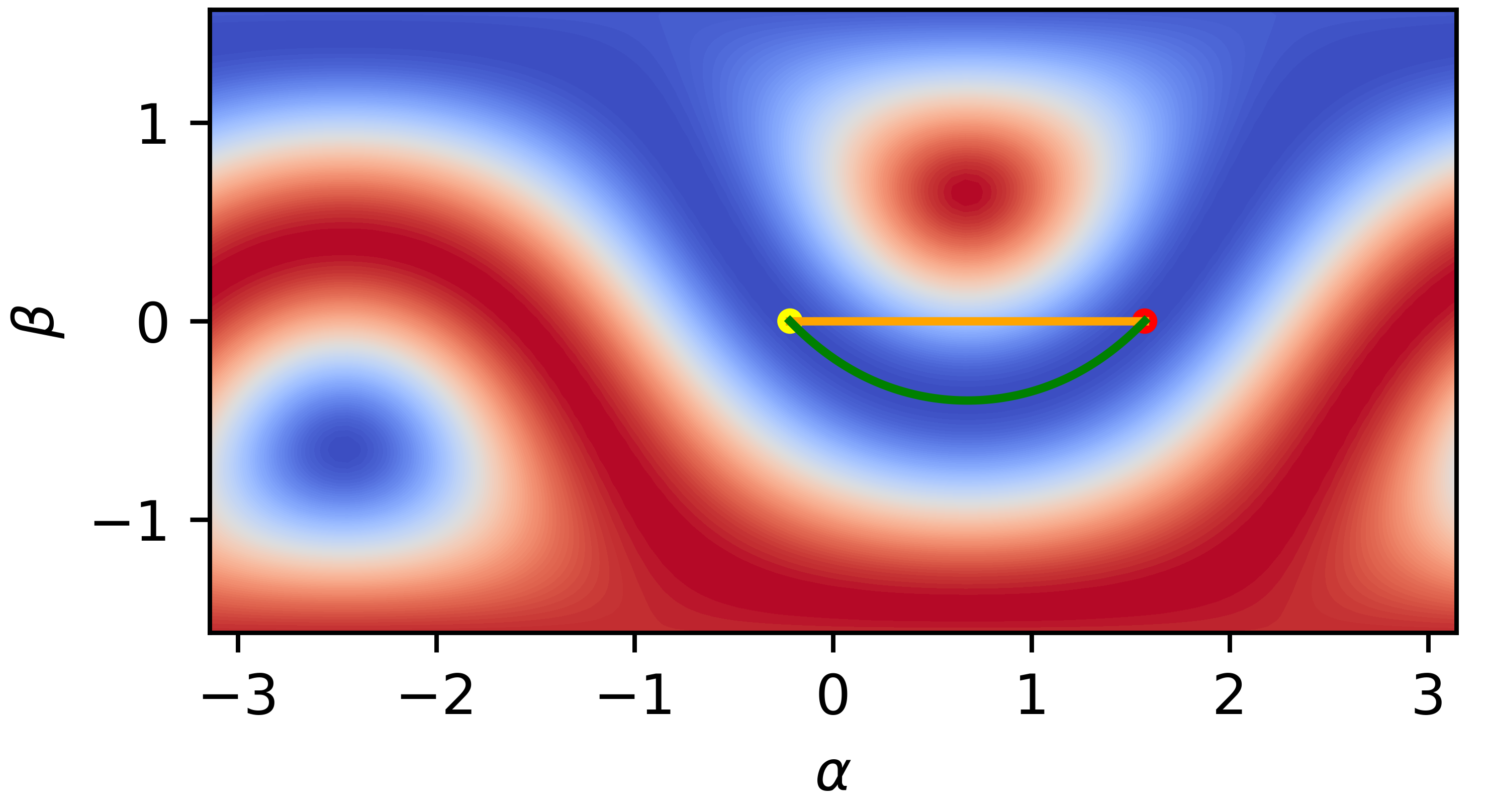} 
\caption{Plot of the optimal paths in terms of $\alpha$ and $\beta$ in the case $\lambda=0$ (orange) and $\lambda=\infty$ (green). The background color represents the value of $\tr(Q^3)$ evaluated at the $Q-$tensor \eqref{Q:alpha-beta}, the uniaxial states ($f(Q)=0$) are indicated in dark blue. The red dot corresponds to $\overline{Q_\infty}$, and the yellow dot \answ{to} $\frac{Q_v}{|Q_v|}$.}
\label{fig:alpha-beta-plane}
\end{center}
\end{figure}

\begin{lemma}\label{lem:lambda:zero:shape-minimizers}
    If $Q$ minimizes $F_0$ with $Q(1)=Q_v$ and $Q(+\infty)=Q_\infty$, then $Q$ is of the form $Q(r)=N(r)\left(\cos(\alpha(r))\Qvb+\sin(\alpha(r))\Qib\right)$, where $\Qib=\sqrt{\frac{3}{2}}Q_\infty$ and
    \begin{equation*}
        \Qvb=\frac{Q_v-\langle Q_v,\Qib\rangle\Qib}{|Q_v-\langle Q_v,\Qib\rangle\Qib|}.
    \end{equation*}
    Moreover, $\alpha$ solves the following ODE,
    \begin{equation*}
        \partial_r(N^2\partial_r\alpha)=-\sqrt{\frac{2}{3}}\cos(\alpha)
    \end{equation*}
    and has boundary values $\alpha(1)=\sin^{-1}(\langle \frac{Q_v}{|Q_v|},\Qib\rangle)$ and $\alpha(+\infty)=\pi/2$, see also Figure \ref{fig:alpha-beta-plane}.
\end{lemma}

\begin{proof}
    For $Q\in Q_\infty+H^1((1,\infty);\sym)$ we use $N=|Q|$ and $\hat{Q}=Q/|Q|$ as seen previously and rewrite the energy $F_0(Q)$ as
    \begin{align*}
        F_0(Q)&=H(N)+K_N(\hat{Q})\\
        &=\int_1^\infty\frac12\abs{\partial_r N}^2\, dr+\int_1^\infty\frac{1}{2}N^2|\partial_r \hat{Q}|^2+g(\hat{Q})\, dr,
    \end{align*}
    where the quantities $H(N)$ and $K_N(\hat{Q})$ are much simpler than in \eqref{H-Kn} since $\lambda=0$. We consider next the equivalent minimization problem,
    \begin{equation*}
        D_0(Q_v) 
        \ = \ 
        \min_{\substack{N\in\sqrt{\frac23}+H^1(1,\infty)\\ N(1) = \sqrt{\frac{2}{3}}}}\left(H(N)+\min_{\hat{Q}}K_N(\hat{Q})\right)
        \, ,
    \end{equation*}
    where the minimum of $K_N(\hat{Q})$ is taken over $\hat{Q}\in \Qib+H^1(1,\infty)$ such that $|\hat{Q}|^2=1$ and $\hat{Q}(1)=\Qvb$. 
    We will begin by addressing the minimization problem for $K_N(\hat{Q})$. 
    We propose that a minimizer will follow the geodesic path on $\S^4$ which connects $\overline{Q_\infty}$ to $\Qvb$. 
    So the minimizer $Q^*$ will be of the form $Q^*(r)=\cos(\alpha(r))\overline{Q_v}+\sin(\alpha(r))\overline{Q_\infty}$. 
    For $Q^*$ to be minimizing amongst maps of this form, $\alpha$ must satisfy the Euler-Lagrange equation associated to
    \begin{equation*}
        K_N(Q^*)=\int_1^\infty\frac{1}{2}N^2|\partial_r\alpha|^2+\sqrt{\frac{2}{3}}(1-\sin(\alpha))\ dr.
    \end{equation*}
    Thus we require $\alpha$ to satisfy,
    \begin{equation}\label{ode:alpha:zero}
        \partial_r(N^2\partial_r\alpha)=-\sqrt{\frac{2}{3}}\cos(\alpha),
    \end{equation}
    with boundary values $\alpha(1)=\alpha_0$ and $\alpha(r)\to\pi/2$ as $r\to\infty$, where $\alpha_0=\sin^{-1}(\langle \frac{Q_v}{|Q_v|},\Qib\rangle)$.
    Next we look at the Euler-Lagrange equation for $K_N$ but subject only to the constraints that $|\hat{Q}|^2=1$ and that $\tr(\hat{Q})=0$. This gives the equation
    \begin{equation}\label{ode:zero}
        -\partial_r(N^2\partial_r\hat{Q})+\hat{Q}_{33}\hat{Q}-N^2|\partial_r\hat{Q}|^2\hat{Q}=Q_\infty.
    \end{equation}
    We show that our candidate minimizer $Q^*$ actually solves the ODE \eqref{ode:zero}. 
    We begin by computing $\partial_r(N^2\partial_r Q^*)$.
    \begin{align*}
        \partial_r(N^2\partial_r Q^*)&=-\partial_r(N^2(\partial_r\alpha)(-\sin(\alpha)\Qvb+\cos(\alpha)\Qib))\\
        &=\partial_r(N^2\partial_r\alpha))(-\sin(\alpha)\Qvb+\cos(\alpha)\Qib)-N^2|\partial_r\alpha|^2(-\cos(\alpha)\Qvb-\sin(\alpha)\Qib)\\
        &=-\sqrt{\frac{2}{3}}\cos(\alpha)(-\sin(\alpha)\Qvb+\cos(\alpha)\Qib)+N^2|\partial_r Q^*|^2Q^*
        \, ,
    \end{align*}
    where we used \eqref{ode:alpha:zero} to get to the final line. Now,
    \begin{align*}
        -\partial_r&(N^2\partial_r Q^*)+Q^*_{33}Q^*-N^2|\partial_r Q^*|^2Q^*\\
        &=\sqrt{\frac{2}{3}}\cos(\alpha)(-\sin(\alpha)\Qvb+\cos(\alpha)\Qib)+\sqrt{\frac{2}{3}}\sin(\alpha)(\cos(\alpha)\Qvb+\sin(\alpha)\Qib)=Q_\infty
        \, .
    \end{align*}
    Therefore to conclude that the minimizer must be of this form, it remains to show that $K_N$ has a unique minimizer. \\

    We note that taking the minimum over $N\in\sqrt{\frac{2}{3}}+H^1(1,\infty)$ is the same as taking the infimum over $N\in\sqrt{\frac{2}{3}}+C_c^\infty(1,\infty)$ due to the density of $C_c^\infty(1,\infty)$ in $H^1(1,\infty)$. 
    Therefore we can treat $N$ as being smooth and compactly supported and by density, we can assume $N>0$. 
    By local existence and uniqueness of geodesics, we know that for $t>0$ small enough $Q^*$ is the unique minimizer of $K_N$.
    The ODE \eqref{ode:zero} can be rewritten as
    \begin{align*}
        -2N(\partial_r N)\partial_r\hat{Q}-N^2\partial_r^2\hat{Q}+\hat{Q}_{33}\hat{Q}-N^2|\partial_r\hat{Q}|^2\hat{Q}=Q_\infty.
    \end{align*}
    and letting $P=\partial_r\hat{Q}$, one can write the second order ODE as the following system of first order ODEs
    \begin{equation*}
        \begin{cases}
            \partial_t\hat{Q} \ = \ P \, , \\
            \partial_t P \ = \ N^{-2}(-2N(\partial_t N)P+\hat{Q}_{33}\hat{Q}-N^2|P|^2\hat{Q}-Q_\infty) \, .
        \end{cases}
    \end{equation*}
    From the uniform Lipschitz continuity of the RHS in $\widehat{Q}$ and $P$, we infer by classical uniqueness for first order ODE systems that the solution $Q^*$ can be uniquely extended for all times.
    
    Thus a minimizer of $K_N(\hat{Q})$ must be of the form $Q^*=\cos(\alpha)\Qvb+\sin(\alpha)\Qib$, so a minimizer of $F_0$ is of the form $Q(r)=N(r)Q^*(r)$.
    \end{proof}

    Now we are ready to prove that the boundary condition associated to $e_\phi$ is minimizing.
    
    \begin{proof}[Proof of Proposition \ref{prop:lambda:zero}]
    We note that we can restrict our attention to $\omega\in\sp\setminus\{\pm e_3\}$ since for $\omega=\pm e_3$, by symmetry of the problem, all admissible $Q_b$ have the same energy.

    So assume that $\omega\in\sp\setminus\{\pm e_3\}$.
    We are using the notation of Lemma \ref{lem:lambda:zero:shape-minimizers} and  claim that $\alpha$ as defined monotonically increases from $\alpha_0$ to $\pi/2$. 
    First note that $\alpha(r)\leq\pi/2$ for all $r\in(1,\infty)$ since otherwise, there would exist a point $r_0$ such that $\alpha(r_0)=\pi/2$ and $\alpha(r_0+\delta)>\pi/2$ for small $\delta$. However if this were true, we could define $\tilde{\alpha}$ such that $\tilde{\alpha}=\alpha$ for $r\leq r_0$ and $\tilde\alpha=\pi/2$ for $r>r_0$. This would result in $\tilde Q=\cos(\tilde\alpha)\Qvb+\sin(\tilde\alpha)\Qib$ having less energy than $Q^*$, but this is not possible since $Q^*$ is minimizing. Thus $\alpha(r)\leq\pi/2$ for all $r\in(1,\infty)$.

    \answ{Assume for the sake of contradiction} that $\alpha$ is not monotone.
    Then, there exists $r_1\in(1,\infty)$ such that $\alpha'(r_1)<0$, but since $\alpha(+\infty)=\pi/2$, there must exist a point $r_2>r_1$ such that $\alpha(r_1)=\alpha(r_2)$ and $\alpha(r)<\alpha(r_1)$ for all $r\in(r_1,r_2)$. We can then define a competitor $\tilde\alpha$ such that $\tilde\alpha=\alpha$ for $r\notin(r_1,r_2)$ and $\tilde\alpha=\alpha(r_1)$ for $r\in(r_1,r_2)$. Then taking $\tilde Q$ as we did previously, only \answ{this time} with a newly defined $\tilde\alpha$, we have lowered the energy, which is not possible as $Q^*$ is a minimizer. Therefore $\alpha$ must monotonically increase from $\alpha_0$ to $\pi/2$.

    Let $\alpha_\phi=\sin^{-1}(\langle\frac{Q_\phi}{|Q_\phi|},\Qib\rangle)$, then $\alpha_\phi\geq\alpha_0$ for all $v\cdot\omega=0$. So by monotonicity of $\alpha$, there exists $r_0\geq 1$ such that $\alpha(r_0)=\alpha_\phi$. 
    Define $\hat\alpha$ by (see Figure \ref{fig:cutting_alpha})
    \begin{equation*}
        \hat{\alpha}(r)=\begin{cases}
            \alpha_\phi,& r\leq r_0,\\
            \alpha(r), & r>r_0.
        \end{cases}
    \end{equation*}

    \begin{figure}[h]
    \begin{center}
    \begin{tikzpicture}[scale=1]
\pgfmathsetmacro{\xmax}{9} 	
\pgfmathsetmacro{\ymax}{5} 	
\pgfmathsetmacro{\t}{0.3} 	
\pgfmathsetmacro{\a}{0.3} 	
\pgfmathsetmacro{\da}{0.26} 	

\draw[->] (0, 0) -- (1.1*\xmax, 0) node[right] {$r$};
\draw[->] (0, 0) -- (0, 1.1*\ymax);

\draw (0,0.1)--(0,-0.1) node[below] {$1$};
\draw (\t*\xmax,0.1)--(\t*\xmax,-0.1) node[below] {$r_0$};
\draw (\xmax,0.1)--(\xmax,-0.1) node[below] {$\infty$};

\draw (0.1,\a*\ymax)--(-0.1,\a*\ymax) node[left] {$\alpha_0$};
\draw (0.1,{(\a+\da)*\ymax})--(-0.1,{(\a+\da)*\ymax}) node[left] {$\alpha_\phi$};
\draw[dotted] (0.1,\ymax)--(-0.1,\ymax) node[left] {$\frac{\pi}{2}$};

\draw[dotted] (\t*\xmax,\ymax)--(\t*\xmax,0); 
\draw[dotted] (\xmax,\ymax)--(\xmax,0); 

\draw[dotted] (0,\ymax)--(\xmax,\ymax); 
\draw[dotted] (0,\a*\ymax)--(\xmax,\a*\ymax); 

\draw[scale=1.0, domain={0}:{\xmax}, samples=100, smooth, variable=\x, line width=1, red!50!black] plot ({\x}, {(\ymax*(\x/\xmax) + \a*\ymax*(1-\x/\xmax))*(1 + 0.1*sin(deg(2*pi*\x/\xmax)))});
\node[red!50!black] at (0.55*\t*\xmax,1.2*\a*\ymax) {$\alpha(\answ{r})$};

\draw[scale=1.0, domain={0}:{\t*\xmax}, samples=100, smooth, variable=\x, line width=1, red] plot ({\x}, {(\a+\da)*\ymax});
\draw[scale=1.0, domain={\t*\xmax}:{\xmax}, dashed, samples=100, smooth, variable=\x, line width=1, red] plot ({\x}, {(\ymax*(\x/\xmax) + \a*\ymax*(1-\x/\xmax))*(1 + 0.1*sin(deg(2*pi*\x/\xmax)))});
\node[red] at (0.55*\t*\xmax,{(\a+\da+0.05)*\ymax}) {$\widehat{\alpha}(\answ{r})$};

\draw[fill=black] (\xmax,\ymax) circle (1.5pt) node[above] {$Q=Q_\infty$};

\end{tikzpicture}
    \caption{Construction of the angle $\widehat{\alpha}$ as a competitor to $\alpha$ in the proof of Proposition \ref{prop:lambda:zero}.}
    \label{fig:cutting_alpha}
    \end{center}
    \end{figure}
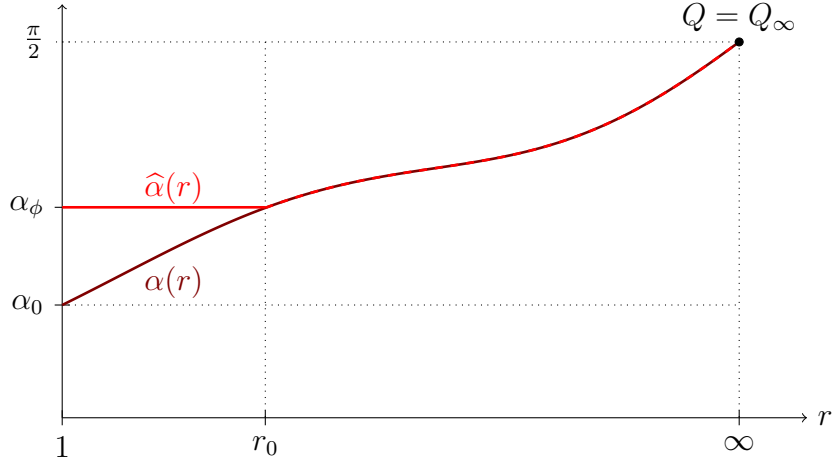

    Let $Q_\phi^*(r)=\cos(\hat{\alpha}(r))\overline{Q_\phi}+\sin(\hat{\alpha}(r))\Qib$, then
    \begin{align*}
        \min_{\hat{Q}(1)=\overline{Q_\phi}}K_N(\hat{Q})
        \ &\leq \
        K_N(Q_\phi^*)
        \ = \
        \int_1^\infty\frac{1}{2}N^2|\partial_r\hat{\alpha}|^2+\sqrt{\frac{2}{3}}(1-\sin(\hat\alpha))\, dr \\
        \ &\leq \
        K_N(Q^*)
        \ = \
        \min_{\hat{Q}(1)=\Qvb}K_N(\hat{Q})
        \, .
    \end{align*}
    It then follows that $D_0(Q_\phi)\leq D_0(Q_v)$ for all $v\cdot\omega=0$.
\end{proof}

\begin{proof}[Proof of Theorem \ref{thm:optim}.]
    The statement of the theorem follows directly from Proposition \ref{prop:min:inf} for $\lambda=\infty$ and Proposition \ref{prop:lambda:zero} for $\lambda=0$.
\end{proof}

Next we consider the case when $\lambda\in (0,\infty)$. 
Here we are still able to prove a specific form for the minimizer along a ray, but we cannot show in full generality that $Q_\phi$ is the minimizing initial condition, see Remark \ref{rem:lambda:finite:problems}.

Note that in the Lyuksyutov regime \cite{KraVi2001,Lyu1978} (i.e.\ for $b=0$ in \eqref{eq:def:f}), the energy $K_N$ for $\lambda>0$ is the same as for $\lambda=0$, so Proposition \ref{prop:lambda:zero} implies that $Q_\phi$ is also the minimizing boundary condition in that case.

\begin{proposition}\label{prop:lambda:finite}
    Let $Q_v=v\otimes v-\frac13 I \in\sym$ for some $v\in\sp$ and let $\lambda\in(0,\infty)$.
    We define $Q_{\mathrm{mix}}=v\otimes e_3+e_3\otimes v-\frac{2}{3}(v\cdot e_3)I$ and $\Qmb$ as
    \begin{align*}
        \Qmb
        \ := \
        \frac{Q_{\mathrm{mix}}-\langle Q_{\mathrm{mix}},\Qbb\rangle-\langle Q_{\mathrm{mix}},\Qib\rangle\Qib}{|Q_{\mathrm{mix}}-\langle Q_{\mathrm{mix}},\Qbb\rangle-\langle Q_{\mathrm{mix}},\Qib\rangle\Qib|}
        \, .
    \end{align*}
    Then minimizers $Q$ of $F_\lambda$ with initial value $Q(1)=Q_v$ are of the form $Q=NQ^*$ where $Q^*$ is given by
    \begin{align}\label{Q:alpha-beta}
        Q^*(r)=\cos(\alpha(r))\cos(\beta(r))\Qbb+\sin(\alpha(r))\cos(\beta(r))\Qib+\sin(\beta(r))\Qmb
        \, ,
    \end{align}
    with
    \begin{equation}\label{bc:alpha-beta}
        \alpha(1)=\sin^{-1}\Big(\Big\langle \frac{Q_v}{|Q_v|},\Qib\Big\rangle\Big),\  \beta(0)=0\quad  \mathrm{and}\quad \alpha(r)\to \pi/2,\ \beta(r)\to 0 \quad\mathrm{as}\quad r\to\infty
        \, .
    \end{equation}
\end{proposition}

In Figure \ref{fig:alpha-beta-plane} we plot the minimizers of $\lambda=0,\infty$ expressed in terms of $\alpha$ and $\beta$ via formula \ref{Q:alpha-beta} and the corresponding value of $\tr((Q^*)^3)$.

\begin{proof}
    Following the same method as in Lemma \ref{lem:lambda:zero:shape-minimizers} we begin with the energy,
    \begin{equation*}
        K_N(Q^*)=
        \int_1^\infty \frac12 N^2(|\partial_r\alpha|^2\cos^2\beta+|\partial_r\beta|^2)+\lambda^2\Big(\frac{2}{9}-N^3\tr((Q^*)^3)\Big)+\sqrt{\frac{2}{3}}(1-\sin\alpha\cos\beta)\, dr
        \, .
    \end{equation*}
    To simplify the notation and since $\tr((Q^*)^3)$ depends only on $\alpha,\beta$ and $v_3$, we define
    \begin{align}
        T(\alpha,\beta,v_3)
        \ &:= \ 
        \tr((Q^*)^3) \notag\\
        \ &= \ \frac{1}{\sqrt6}\sin^3(\alpha)\cos^3(\beta)+\frac{1}{2\sqrt{2}}\left(\frac{(v_3^2-1)(9v_3^2-1)}{(1-v_3^2)^{1/2}(3v_3^2+1)^{3/2}}\right)\cos(\alpha)\sin^2(\beta)\cos(\beta)\notag\\
        &+\frac{1}{\sqrt{6}}\left(\frac{v_3^2-\frac{1}{3}}{v_3^2+\frac{1}{3}}\right)\cos^2(\alpha)\sin(\alpha)\cos^3(\beta)+\frac{1}{\sqrt{6}}\left(\frac{1-v_3^2}{(v_3^2+\frac{1}{3})^{3/2}}\right)\sin^3(\beta)\notag\\
        &+\sqrt{\frac{2}{3}}\left(\frac{v_3^2-\frac{1}{3}}{(v_3^2+\frac{1}{3})^{3/2}}\right)\cos^2(\alpha)\sin(\beta)\cos^2(\beta)\label{trace}\\
        &+\frac{1}{2\sqrt{6}}\left(\frac{1-9v_3^2}{3v_3^2+1}\right)\sin(\alpha)\sin^2(\beta)\cos(\beta)+\sqrt{\frac{2}{3}}\left(\frac{v_3^2(1-v_3^2)^{1/2}}{(v_3^2+\frac{1}{3})^{3/2}}\right)\cos^3(\alpha)\cos^3(\beta)\notag\\
        &-\left(\frac{v_3\sqrt{6-6v_3^2}}{6v_3^2+2}\right)\sin(\alpha)\cos(\alpha)\sin(\beta)\cos^2(\beta)\notag
        \, .
    \end{align}
    Then we obtain the Euler-Lagrange equations,
    \begin{equation}\label{el:alpha}
        \partial_r(N^2(\partial_r\alpha)\cos^2\beta)=-\lambda^2N^3\frac{\partial T}{\partial \alpha}-\sqrt{\frac{2}{3}}\cos\alpha\cos\beta
        \, ,
    \end{equation}
    \begin{equation}\label{el:beta}
        \partial_r(N^2\partial_r\beta)=-N^2|\partial_r\alpha|^2\sin\beta\cos\beta-\lambda^2N^3\frac{\partial T}{\partial \beta}+\sqrt{\frac23}\sin\alpha\sin\beta
        \, .
    \end{equation}
    So by minimizing amongst maps of the same form as $Q^*$, we can see that the minimizer must have $\alpha$ and $\beta$ which solve equations \eqref{el:alpha} and \eqref{el:beta} subject to the boundary values \eqref{bc:alpha-beta}. 
    From the energy $K_N(\hat{Q})$, and due to the constraints $|\hat{Q}|^2=1$ and $\tr(\hat{Q})=0$, we get the following Euler-Lagrange equation including Lagrange multipliers,
    \begin{equation}\label{el:lambda-finite}
        -\partial_r(N^2\partial_r\hat{Q}) -3\lambda^2N^3\Big(\hat{Q}^2-\tr(\hat{Q}^3)-\frac{1}{3}I\Big)+\hat{Q}_{33}\hat{Q}-N^2|\partial_r \hat{Q}|^2\hat{Q}=Q_\infty
        \, .
    \end{equation}
    To prove that $Q^*$ is the unique minimizer of $K_N$, we  must first demonstrate that it satisfies \eqref{el:lambda-finite}. 
    This can be done through a direct but tedious computation. First, we can rewrite \eqref{el:alpha} as,
    \begin{equation}\label{el:alpha:new}
        \partial_r(N^2(\partial_r\alpha)\cos\beta)=N^2(\partial_r\alpha)(\partial_r\beta)\sin\beta-\lambda^2N^3(\cos\beta)^{-1}\frac{\partial T}{\partial \alpha}-\sqrt{\frac{2}{3}}\cos\alpha
        \, ,
    \end{equation}
    where we note that since every term in $\frac{\partial T}{\partial\alpha}$ has a factor of $\cos\beta$, we can divide by $\cos\beta$. 
    We then evaluate $-\partial_r(N^2\partial_r Q^*)$, using equations \eqref{el:beta} and \eqref{el:alpha:new} to obtain
    \begin{align*}
        -\partial_r(N^2\partial_r Q^*)&=-Q^*_{33}Q^*+N^2|\partial_r Q^*|^2Q^*-\lambda^2N^3\Big((\cos\beta)^{-1}\frac{\partial T}{\partial \alpha}(\sin(\alpha)\Qvb-\cos(\alpha)\Qib)\Big)\\
        &\quad-\lambda^2N^3\Big(\frac{\partial T}{\partial \beta}\big(\sin(\beta)(\cos(\alpha)\Qbb+\sin(\alpha)\Qib)-\cos(\beta)\Qmb\big)\Big)+Q_\infty
        \, .
    \end{align*}
    This is done by expanding and collecting terms to simplify the expression. 
    From here it suffices to show that
    \begin{align}
        (\cos\beta)^{-1}&\frac{\partial T}{\partial \alpha}(\sin(\alpha)\Qvb+\cos(\alpha)\Qib)\label{T-computation}\\
        -&\frac{\partial T}{\partial \beta}\big(\sin(\beta)(\cos(\alpha)\Qbb+\sin(\alpha)\Qib)-\cos(\beta)\Qmb\big)
        +3(Q^*)^2-3T Q^*-I=0
        \, .\notag
    \end{align}
    As $T$ is given by \eqref{trace}, the calculations of its derivatives are quite involved. Nevertheless, with the help of the software sympy, we can show  that \eqref{T-computation} holds and thus $Q^*$ solves \eqref{el:lambda-finite}.
    
    To show uniqueness of the solution to the ODE, the exact same argument as in the $\lambda=0$ case can be applied, only with a slightly more complicated system of ODEs,
    \begin{equation*}
        \begin{cases}
            \partial_r\hat{Q}=P,\\
            \partial_r P=N^{-2}\big(-2NP(\partial_r N)-3\lambda^2N^3(\hat{Q}^2-\tr(\hat{Q}^3)-\frac{1}{3}I)+\hat{Q}_{33}\hat{Q}-N^2|P|^2\hat{Q}-Q_\infty\big)\, .
        \end{cases}
    \end{equation*}
    Thus $K_N(\hat{Q})$ is minimized by $Q^*$ and we can subsequently conclude that a minimizer of $F_\lambda$ is of the form $Q=NQ^*$.
\end{proof}
\begin{remark}\label{rem:lambda:finite:problems}
    The equations for $\alpha$ and $\beta$ in Proposition \ref{prop:lambda:finite} are coupled and considerably more difficult to analyze than in the case $\lambda=0$. 
    For example, our previous argument to obtain monotonicity for $\alpha$ does not apply anymore.
    The construction of a competing $\alpha$ is no longer guaranteed to decrease the energy due to the lack of monotonicity of $\tr(Q^3)$ in $\alpha$ (see \eqref{trace} and Figure \ref{fig:alpha-beta-plane}).
    In addition, the energy for $\lambda\in (0,\infty)$ also depends on $v_3$ via $\tr(Q^3)$ and thus one would also need monotonicity in $v_3$ as we change the boundary condition alongside $\alpha$. Nevertheless we conjecture that $Q_\phi$ is the minimizer of $D_\lambda$ also in the case $\lambda\in (0,\infty)$ and the trajectory of the optimal profile $Q$ remains in between those of the optimal $Q$'s for $\lambda=0$ and $\lambda=\infty$, see again Figure \ref{fig:alpha-beta-plane}.
\end{remark}

\section{Upper Bound}
\label{sec:upper}

In this section we prove the upper bound in Theorem \ref{thm:upper} by constructing a sequence of competitors for the cases $\lambda=\infty$ and $\lambda\in [0,\infty)$ separately. 
Theorem \ref{thm:upper} thus follows immediately from Proposition \ref{prop:up:inf} and Proposition \ref{prop:up:fin}.

\subsection{The Case $\lambda=\infty$}
\label{ss:upI}

We begin with the case where $\lambda=\infty$, and we wish to prove the following proposition.
\begin{proposition}\label{prop:up:inf}
    If $\Qm$ minimizes $\exe$ with boundary condition $\Qm=\answ{Q_\phi}$ on $\mathbb S^2$ and
        \[\frac{\eta}{\xi}\to\infty,\quad\text{as}\quad \xi\dth{,\eta}\to0,\]
        then,
        \[\limsup_{\xi\dth{,\eta}\to0}\eta \exe(\Qm)\leq\int_{\mathbb S^2}D_\infty(\answ{Q_\phi}(\omega))\ \answ{d\H^2(\omega)}. \]
\end{proposition}
To do this, we will construct a sequence maps which are uniaxial a.e., equivariant and symmetric with respect to $\{x_3=0\}$.

\begin{proof}[Proof of Proposition \ref{prop:up:inf}]
    To simplify the construction, we will construct the recovery sequence on
    $\Om_0^+:=\{(r,\phi):1<r<\infty,\ 0<\phi<\frac{\pi}{2}\}$ 
    and then extend it to $\Om$ via a rotation, yielding an equivariant map. For the construction we partition $\Om_0^+$ into smaller regions $\Omega_k$, $k=1,...,4$ (see Figure \ref{fig:infty}).
    \answ{In each region, we will construct a competitor and show that the only contribution to the energy in the limit $\xi,\eta\to 0$ comes from the region $\Omega_1$.}
    Since the recovery sequence will consist of uniaxial maps, it is enough to define a continuous and piecewise smooth map $\hat{n}:\Om_0^+\to\S^2$ and then take
    \begin{align*}
        \hat{Q}_{\xi\dth{,\eta}}(r,\phi)=\begin{cases}
            \hat{n}(r,\phi)\otimes \hat{n}(r,\phi)-\frac{1}{3}I,& 0<\phi<\frac{\pi}{2},\\
            (T\circ\hat{n}(r,\pi-\phi))\otimes (T\circ\hat{n}(r,\pi-\phi))-\frac{1}{3}I,& \frac{\pi}{2}<\phi<\pi,
        \end{cases}
    \end{align*}
    where $T:\sp\to\sp$ is the reflection $T(x_1,x_2,x_3)=(x_1,x_2,-x_3)$.
    The competitor sequence is then $\overline{\Qm}(r,\theta,\phi)=R_\theta^T\hat{Q}_{\xi\dth{,\eta}}(r,\phi) R_\theta$. Note that we will use the notation, $\Omega_k':=\{(r,\theta,\phi):(r,\phi)\in\Omega_k,\ 0\leq\theta<2\pi\}$ to allow us to define the maps on $\Om_0^+$, but consider the energy on the entire upper space $\Om^+$. \\

    \begin{figure}[h]
    \begin{center}
    \begin{tikzpicture}[scale=1.5]

    \draw (0,1) arc (90:0:1);
    \draw (0,1) -- (0,3);
    \draw (1,0) -- (3,0);
    \draw (45:1)-- ++(45:3);
    \draw (3*45/2:1)-- ++(3*45/2:0.5);
    \draw (90:1.5) arc (90:3*45/2:1.5);
    \draw (90:2) arc (90:45:2);
    \draw[<->] (0,0.9) arc (90:45:.9);
    \draw (70:0.8) node [below] {$2\eta$};
    \draw[<->] (40:1.05) -- ++(45:0.9);
    \draw (41:1.4) node [below right] {$2\eta$};
    \draw (2.5,1) node {$\Omega_1$};
    \draw (1,2.5) node {$\Omega_2$};
    \draw (0.25,1.2) node {$\Omega_3$};
    \draw (0.9,1.4) node {$\Omega_4$};
    \end{tikzpicture}
    \end{center}
    \caption{Subdivision of $\Om_0^+$ for the proof of Proposition \ref{prop:up:inf}}\label{fig:infty}
    \end{figure}
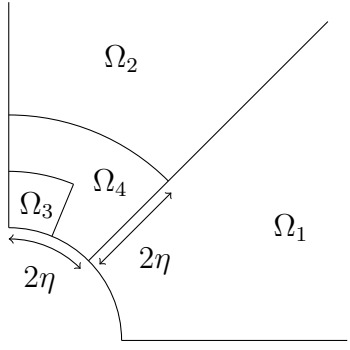

    \noindent \underline{\textit{Energy in $\Omega_1$}}: We begin by defining \[\Omega_1:=\left\{(r,\phi):1<r<\infty,\ 2\eta<\phi<\frac{\pi}{2}\right\},\]
    and take $\hat{n}(r,\phi)=n\left(\frac{r-1}{\eta},\frac{\pi}{2}-\phi\right)$ where $n$ is given as in Lemma \ref{def:n_min}. Note that the use of $\frac{\pi}{2}-\phi$ is due to the fact that $\phi$ is representing the point on the sphere instead of the angle between the tangent vector and $e_3$. \dth{Expressing the energy in spherical coordinates $(r,\theta,\phi)$ and using equivariance}
    \begin{align*}
        \eta \exe(\overline{\Qm}\dth{;}\Omega_1')=2\pi\eta\int_{2\eta}^\frac{\pi}{2}\int_1^\infty r^2\sin\phi\left|\frac{\partial \hat{n}}{\partial r}\right|^2+\sin\phi\left|\frac{\partial \hat{n}}{\partial \phi}\right|^2+\frac{|\hat{n}_1|^2}{\sin\phi}+\frac{r^2\sin\phi}{\eta^2}g(\hat{n})\ drd\phi
        \, .
    \end{align*}
    Then by change of variables $r=1+\eta t$,
    \begin{multline*}
        \eta \exe(\overline{\Qm}\dth{;}\Omega_1')=2\pi\int_{2\eta}^\frac{\pi}{2}\int_0^\infty (1+\eta t)^2\sin\phi\left|\frac{\partial n}{\partial t}\right|^2+\eta^2\sin\phi\left|\frac{\partial n}{\partial \phi}\right|^2+\frac{\eta^2|n_1|^2}{\sin\phi}\\
        +(1+\eta t)^2\sin\phi\  g(n)\ dtd\phi.
    \end{multline*}
    We divide this integral into two parts: 
    \begin{align}\label{split}
        \eta \exe(\Qcon\dth{;}\Om_1')=2\pi\int_{2\eta}^\frac{\pi}{2}\left(\int_0^\infty \abs{\frac{\partial n}{\partial t}}^2+g(n)\ dt\right)\sin\phi\ d\phi+\mathcal{R}_\dth{\eta}
    \end{align}
    where,
    \begin{multline*}
        \mathcal{R}_\dth{\eta}=2\pi\int_{2\eta}^\frac{\pi}{2}\dth{\int_0^\infty\bigg[}(2\eta t+\eta^2t^2)\sin\phi\left|\frac{\partial n}{\partial t}\right|^2+\eta^2\sin\phi\left|\frac{\partial n}{\partial \phi}\right|^2+\frac{\eta^2|n_1|^2}{\sin\phi}\\
        +(2\eta t+\eta^2t^2)\sin\phi\  g(n)\dth{\bigg]}\ dtd\phi.
    \end{multline*}
    We identify the first term in \eqref{split} as \answ{$D_\infty(Q_\phi)$ using Lemma \ref{def:n_min}}, while the second term $\mathcal{R}_\dth{\eta}$ vanishes in the limit. For this second part,
    we note that since $2\eta<\phi<\frac{\pi}{2}$, we still have $0<\frac{\pi}{2}-\phi<\frac{\pi}{2}$ so we can use the bounds from \eqref{n_min:bound} to get,
    \begin{align*}
        \mathcal{R}_\dth{\eta}\leq 2\pi\int_{2\eta}^\frac{\pi}{2}\int_0^\infty\dth{\bigg[}(2\eta t+\eta^2t^2)Ce^{-\kappa t}+C\eta^2e^{-\kappa t}+\frac{C\eta^2e^{-\kappa t}}{\sin 2\eta}+(2\eta t+\eta^2t^2)\sqrt{\frac{3}{2}}Ce^{-\kappa t}\dth{\bigg]}\ dtd\phi.
    \end{align*}
    Then by straightforward computations, we get,
    \begin{align*}
        \mathcal{R}_\dth{\eta}\leq2\pi\left(\frac{\pi}{2}-2\eta\right)\left(C\eta+C\eta^2+\frac{C\eta^2}{\sin2\eta}\right)
    \end{align*}
    for some constants $C>0$, which tends to $0$ as $\eta\to 0$. 
    \answ{Hence restating \eqref{split} using $D_\infty(Q_\phi)$ and changing the integration variable from $\phi$ to $\varphi$ so as not to be confused with the notation $Q_\phi$,} it follows that
    \begin{align*}
        \eta \exe(\overline{\Qm}\dth{;}\Omega_1')&=2\pi\int_{2\eta}^{\frac{\pi}{2}}D_\infty(\answ{Q_\phi}(0,\varphi)) \sin \varphi\, d \varphi+\mathcal{R}_\dth{\eta}\leq\int_{\S^2_+}D_\infty(\answ{Q_\phi}(\omega))\ \answ{d\H^2(\omega)}+\mathcal{R}_\dth{\eta},
    \end{align*}
    therefore when taking $\xi\dth{,\eta}\to0$, we are left with
    \begin{align*}
        \limsup_{\xi\dth{,\eta}\to0}\eta \exe(\overline{\Qm}\dth{;}\Om_1')\leq\int_{\S^2_+}D_\infty(\answ{Q_\phi}(\omega))\ \answ{d\H^2(\omega)}.
    \end{align*}
    
    \noindent \underline{\textit{Energy in $\Omega_2$}}: Next consider the region $\Om_2$ (see also Figure \ref{fig:infty}) defined to be\[\Om_2:=\{(r,\phi):1+2\eta<r<\infty,\ 0<\phi<2\eta\}.\]
    Because of equivariance, we want to take $\hat{n}=e_3$ on the $x_3-$axis. Furthermore, at $\phi=2\eta$ we need $\hat{n}$ to be consistent with the construction from the previous region $\Om_1$. To preserve uniaxiality, we define $\hat{n}$ by an interpolation of its angle with $e_3$. On $\phi=2\eta$ we take $\hat{n}$ as defined in $\Om_1$ and at $\phi=0$, we take $\hat{n}=e_3$. 
    Let $\Phi$ denote the angle between $e_3$ and our construction $\hat{n}$ on $\Om_2$, then we define $\Phi$ as follows,
    \begin{equation*}
        \Phi(r,\phi)=\frac{\phi}{2\eta}\cos^{-1}(\hat{n}_3(r,2\eta)),
    \end{equation*}
    where $\hat{n}_3(r,2\eta)$ is defined \dth{by continuous extension onto $\partial \Om_1$}. So we let $\hat{n}=(-\sin\Phi,0,\cos\Phi)$ on $\Om_2$. We note that
    \[\left|\frac{\partial \hat{n}}{\partial r}\right|^2=\left|\frac{\partial \Phi}{\partial r}\right|^2\quad\text{and}\quad \abs{\frac{\partial \hat{n}}{\partial \phi}}^2=\left|\frac{\partial \Phi}{\partial \phi}\right|^2,\]
    so we can write the energy in this region as
    \begin{align*}
        \eta \exe(\Qcon\dth{;}\Omega_2')=2\pi\eta\int_0^{2\eta}\int_{1+2\eta}^\infty r^2\sin\phi\abs{\frac{\partial \Phi}{\partial r}}^2+\sin\phi\abs{\frac{\partial \Phi}{\partial \phi}}^2+\frac{\sin^2\Phi}{\sin\phi}+\frac{r^2\sin\phi}{\eta^2}\sin^2\Phi \ drd\phi.
    \end{align*}
    By the change of variables, $r=1+\eta t$, and by letting $\hat{\Phi}(t,\phi)=\Phi(1+\eta t,\phi)$, we have
    \begin{align*}
        \eta \exe(\Qcon\dth{;}\Omega_2')\leq2\pi\int_0^{2\eta}\int_{2}^\infty (1+\eta t)^2\abs{\frac{\partial \hat{\Phi}}{\partial t}}^2+\eta^2\abs{\frac{\partial \hat{\Phi}}{\partial \phi}}^2+\frac{\eta^2\hat{\Phi}^2}{\sin\phi}+(1+\eta t)^2\hat{\Phi}^2\sin\phi \ dtd\phi.
    \end{align*}
    We note that the $t-$derivative term is bounded by $Ce^{-\kappa t}$ since,
    \[\abs{\frac{\partial \hat{\Phi}}{\partial t}}^2\leq\abs{\frac{\partial}{\partial t}\hat{\Phi}(t,2\eta)}^2=\abs{\frac{\partial n}{\partial t}(t,2\eta)}^2\leq Ce^{-\kappa t}.\]
    by \eqref{n_min:bound} of Lemma \ref{def:n_min}. We get a similar bound on the $\phi-$derivative,
    \begin{align*}
        \abs{\frac{\partial\hat{\Phi}}{\partial \phi}}^2=\frac{1}{4\eta^2}\abs{\cos^{-1}(n_3(t,2\eta))}^2=\frac{1}{4\eta^2}\abs{\hat{\Phi}(t,2\eta)}^2
    \end{align*}
    and using that $\alpha\leq 2\sin\alpha$ for $0<\alpha<\frac{\pi}{2}$ we get,
    \begin{align}
        \abs{\frac{\partial\hat{\Phi}}{\partial \phi}}^2\leq\frac{1}{4\eta^2}\abs{2\sin(\hat{\Phi}(t,2\eta))}^2=\frac{1}{\eta^2}|n_1(t,2\eta)|^2\leq \frac{Ce^{-\kappa t}}{\eta^2}.\label{bd:om2:phi}
    \end{align}
    Using
    \begin{equation*}
        \frac{\eta^2\hat{\Phi}^2}{\sin\phi}=\frac{\phi^2(\cos^{-1}(n_3(t,2\eta))^2}{4\sin\phi}\leq\frac{\phi\left(\hat{\Phi}(t,2\eta)\right)^2}{4}\leq C|n_1(t,2\eta)|^2,
    \end{equation*}
    we have the following bound on the energy,
    \begin{align*}
        \eta \exe(\Qcon\dth{;}\Om_2')&\leq2\pi\int_0^{2\eta}\int_2^\infty Ce^{-\kappa t}((1+\eta t)^2+1)+C|n_1(t,2\eta)|^2+(1+\eta t^2)|n_1(t,2\eta)|^2\ dtd\phi\\
        &\leq2\pi\int_0^{2\eta}\int_2^\infty Ce^{-\kappa t}((1+\eta t)^2+1)\ dtd\phi\, \leq\,  4\pi\eta C,
    \end{align*}
    for some constant $C>0$. Therefore the energy contribution from $\Om_2$ is negligible in the limit $\xi\dth{,\eta}\to0$.
    
    \noindent \underline{\textit{Energy in $\Omega_3$}}: In this region, we must define $\Qcon$ to have a point singularity at the pole of the sphere, due to the discontinuity in the boundary condition at this point, however the energy will still be bounded since point singularities in 3D have finite energy. With the additional $\eta-$prefactor in our energy, the energy contribution from $\Om_3$ will therefore be negligible as well. 
    Let 
    \[\Omega_3:=\{(r,\phi):1<r<1+\eta,\ 0<\phi<\eta\}.\]
    On a flat domain, it would be possible to define $\hat{n}$ to be the ``standard'' point singularity at $(0,0)$,
    \begin{equation}
        m(s,\tau):=\left(\frac{-\tau}{\sqrt{\tau^2+s^2}},0,\frac{s}{\sqrt{\tau^2+s^2}}\right).\label{sing1}
    \end{equation}
    Because our domain is curved, we have to slightly adapt this profile in order to match the boundary conditions. 
    Thus,
    \begin{equation}
        \hat{n}(r,\phi)=R_\phi m\left(\frac{r-1}{\eta},\frac{\phi}{\eta}\right),\quad\text{where}\quad R_\phi=\begin{pmatrix}
            \cos\phi&0&\sin\phi\\
            0&1&0\\
            -\sin\phi&0&\cos\phi
        \end{pmatrix}.\label{sing2}
    \end{equation}
    so we have $r=1+\eta s$ and $\phi=\eta\tau$, thus if $(r,\phi)\in \Om_3$, then $(s,\tau)\in(0,1)\times(0,1)$. We first consider the $r-$derivative term,
    \begin{equation}
        \abs{\frac{\partial \hat{n}}{\partial r}}^2=\frac{1}{\eta^2}\abs{\frac{\partial m}{\partial s}}^2=\frac{1}{\eta^2}\left(\frac{\tau^2}{|(s,\tau)|^4}\right)\label{sing:r}
        \, ,
    \end{equation}
    by a direct computation. Here $|(s,\tau)|$ denotes the standard Euclidean 2-norm. The derivative in $\phi$ is more complicated, due to the rotation in the definition of $\hat{n}$. But again by a direct computation we can show that,
    \begin{equation*}
        \abs{\frac{\partial \hat{n}}{\partial \phi}}^2=|m|^2+\frac{1}{\eta^2}\abs{\frac{\partial m}{\partial \tau}}+\frac{2}{\eta}\left(m_3\frac{\partial m_1}{\partial \tau}-m_1\frac{\partial m_3}{\partial \tau}\right)\leq 1+\frac{1}{\eta^2}\abs{\frac{\partial m}{\partial \tau}}^2.
    \end{equation*}
    After computing the $\tau-$derivative explicitly, we have
    \begin{align}
        \abs{\frac{\partial \hat{n}}{\partial \phi}}^2\leq1+\frac{1}{\eta^2}\left(\frac{s^2}{|(s,\tau)|^4}\right).\label{sing:phi}
    \end{align}
    Next we have $\hat{n}_1=m_1\cos\phi+m_3\sin\phi$ and we notice that since $|m|=1$, it follows that $|\hat{n}_1|^2$ and $g(\Qcon)$ are bounded by a constant. However in order to get the required estimate for the $\theta-$derivative, we need a more exact bound on $|\hat{n}_1|$. By definition of $\hat{n}_1$ and since \dth{$\frac12\phi\leq \sin\phi$,}
    \begin{align*}
        \frac{|\hat{n}_1|^2}{\sin\phi}\leq\frac{2(m_1\cos\phi+m_3\sin\phi)^2}{\phi}=\frac{2(m_1^2\cos^2\phi+2m_1m_3\sin\phi\cos\phi+m_3^2\sin^2\phi)}{\phi}
        \, .
    \end{align*}
    \dth{Then using $\sin\phi\leq\phi$ and} substituting $\phi=\eta\tau$, it holds,
    \begin{equation*}
        \frac{|\hat{n}_1|^2}{\sin\phi}\leq\frac{2m_1^2}{\eta\tau}+6=\frac{C}{\eta}\left(\frac{\tau}{|(s,\tau)|^2}\right)+6
        \, .
    \end{equation*}
    for some constant $C>0$. 
    Now we are ready to compute the energy in this region. Altogether, we observe that the energy in $\Om_3'$ is estimated as follows.
    \begin{align}
        \eta \exe(\Qcon\dth{;}\Om_3')&=2\pi\eta\int_0^\eta\int_1^{1+\eta}r^2\sin\phi\abs{\frac{\partial \hat{n}}{\partial r}}^2+\sin\phi\abs{\frac{\partial \hat{n}}{\partial \phi}}^2+\frac{|\hat{n}_1|^2}{\sin\phi}+\frac{r^2\sin\phi}{\eta^2}g(\Qcon)\ drd\phi\notag\\
        &\leq 2\pi\eta\int_0^1\int_0^1 C\sin(\eta\tau)\left(\frac{\tau^2}{|(s,\tau)|^4}+\frac{s^2}{|(s,\tau)|^4}\right)+\eta^2+\frac{C\eta\tau}{|(s,\tau)|^2}+C\eta\ dsd\tau\notag\\
        &\leq 2\pi\eta\int_0^1\int_0^1C\eta\left(\frac{\tau^3}{|(s,\tau)|^4}+\frac{s^2\tau}{|(s,\tau)|^4}+\frac{\tau}{|(s,\tau)|^2}\right)+C\eta\ dsd\tau\notag\\
        &\leq 2\pi C\eta^2
        \, .\label{sing:en}
    \end{align}
    Taking the limit $\xi\dth{,\eta}\to0$, the energy contribution of $\Omega_3'$ vanishes.
    
    \noindent \underline{\textit{Energy in $\Omega_4$}}: Finally we consider $\Om_4$ where we will define $\hat{n}$ by a Lipschitz extension. Let \[\Omega_4:=\{(r,\phi):1<r<1+2\eta,\ 0<\phi<2\eta\}\setminus \overline{\Omega_3},\]
    then we will first define a function $\Phi$ on $\Omega_4$ and the let $\hat{n}=(\sin\Phi,0,\cos\Phi)$. The boundary $\partial\Omega_4$ can be split into 6 pieces as shown in Figure \ref{om_4}.
    
    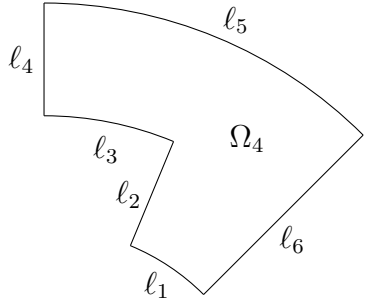
\begin{figure}[h]
    \begin{center}
    \begin{tikzpicture}[scale=1.5]

    \draw (3*45/2:2) arc (3*45/2:45:2);
    \draw (0,3) -- (0,4);
    \draw (45:2)--(45:4);
    \draw (3*45/2:2)--(3*45/2:3);
    \draw (90:3) arc (90:3*45/2:3);
    \draw (90:4) arc (90:45:4);
    \draw (1,1.5) node {$\ell_1$};
    \draw (0.75,2.3) node {$\ell_2$};
    \draw (0.55,2.7) node {$\ell_3$};
    \draw (-0.2,3.5) node {$\ell_4$};
    \draw (1.7,3.85) node {$\ell_5$};
    \draw (2.2,1.9) node {$\ell_6$};
    \draw (1.8,2.8) node {$\Omega_4$};
    \end{tikzpicture}
    \end{center}
    \caption{Boundary components of $\Om_4$ for Proposition \ref{prop:up:inf}}
    \label{om_4}
    \end{figure}
    
    On $\ell_1$ we choose $\Phi$ so that $\hat{n}=e_\phi$ as required by the boundary condition, then on $\ell_2$ and $\ell_3$, we define $\hat{n}$ so that it is consistent with the construction on $\Om_3$. 
    On $\ell_4$ we let $\Phi=0$ so that $\hat{n}=e_3$ and on $\ell_5$ and $\ell_6$ we define $\hat{n}$ to be consistent with $\Om_2$ and $\Om_1$ respectively. 
    We will then show that $\Phi$, being the angle between $\hat{n}$ and $e_3$, is Lipschitz on $\partial \Omega_4$ with Lipschitz constant proportional to $\eta^{-1}$ so we can extend $\Phi$ to all of $\Omega_4$ with the same Lipschitz constant using the Kirszbraun Theorem, see e.g.\ \cite[p.21]{lip}. 
    To show $\Phi$ has this Lipschitz constant, we consider the derivatives along each component of the boundary.

    First on $\ell_1$, $\answ{Q_\phi}$ is Lipschitz with constant $2$, so $\Phi$ is also Lipschitz with a constant that has no dependence on $\eta$. 
    Next on $\ell_2$ and $\ell_3$, we consider the derivatives,
    \begin{equation*}
        \abs{\frac{\partial \hat{n}}{\partial r}}^2=\frac{C}{\eta^2}\left(\frac{\tau^2}{|(s,\tau)|^4}\right)\leq \frac{C}{\eta^2}\quad\text{and}\quad\abs{\frac{\partial \hat{n}}{\partial \phi}}^2=1+\frac{1}{\eta^2}\left(\frac{s^2}{|(s,\tau)|^4}\right)\leq \frac{C}{\eta^2}
        \, ,
    \end{equation*}
    using \eqref{sing:r} and \eqref{sing:phi} since $\tau=1$ and $s=1$ on $\ell_2$ and $\ell_3$ respectively. 
    On $\ell_4$, $\Phi$ is constant, so it has Lipschitz constant 0. For $\ell_5$, we have
    \begin{equation*}
        \abs{\frac{\partial \hat{n}}{\partial \phi}}^2\leq\frac{Ce^{-\kappa t}}{\eta^2}\leq\frac{C}{\eta^2}
        \, ,
    \end{equation*}
    by \eqref{bd:om2:phi} and using that $t=2$ here. Finally, on $\ell_6$, 
    \begin{equation*}
        \abs{\frac{\partial \hat{n}}{\partial r}}^2=\frac{1}{\eta^2}\abs{\frac{\partial n}{\partial t}}^2\leq \frac{Ce^{-\kappa t}}{\eta^2}\leq \frac{C}{\eta^2}
        \, ,
    \end{equation*}
    since $0<t<2$ on this region. So $\Phi$ is Lipschitz with constant $C\eta^{-1}$ on $\partial\Omega_4$ and evoking Kirszbraun Theorem there exists a Lipschitz extension $\Phi$ to all of $\Omega_4$ with the same Lipschitz constant. Thus each derivative of $\Phi$ on $\Omega_4$ is bounded by $C\eta^{-1}$. Using the Lipschitz constant, we get that $|\Phi(\phi)-\Phi(0)|\leq C\eta^{-1}|\phi|$ and therefore,
    \begin{align*}
        \eta \exe(\overline{\Qm}\dth{;}\Omega_4')&\leq2\pi\eta\int_0^{2\eta}\int_\gamma^{1+2\eta}\frac{Cr^2\sin\phi}{\eta^2}+\frac{C\sin\phi}{\eta^2}+\frac{2\phi}{\eta^2}+\sqrt{\frac{3}{2}}\frac{r^2\phi^3}{\eta^2}\ drd\phi\\
        &\leq2\pi\eta\int_0^{2\eta}\int_\gamma^{1+2\eta}\frac{C(1+2\eta)^2}{\eta}+\frac{C}{\eta}+\frac{4}{\eta}+C(1+2\eta)^2\eta\ drd\phi\\
        &\leq2\pi\eta(4\eta^2)\left[\frac{C(1+2\eta)^2}{\eta}+\frac{C}{\eta}+\frac{4}{\eta}+C(1+2\eta)^2\eta\right]
        \, .
    \end{align*}
    Thus when we take $\xi\dth{,\eta}\to0$, the whole energy will vanish on this region.
    
    \noindent\underline{\textit{Energy in $\Om$}}: \answ{Having analyzed} the energy on each region, we are ready to put the regions together and get a bound on the energy in all of $\Om$. 
    We have that
    \begin{align*}
        \eta \exe(\Qm)&\leq\eta \exe(\Qcon)=2\eta \exe(\Qcon\dth{;}\Om_+)=2\sum_{k=1}^4\eta \exe(\Qcon\dth{;}\Om_k')
        \, ,
    \end{align*}
    by symmetry of the construction. So applying our results from each region, we can see that
    \begin{align*}
        \limsup_{\xi\dth{,\eta}\to0}\eta \exe(\Qm)&\leq2\sum_{k=1}^4\limsup_{\xi\dth{,\eta}\to0}\eta \exe(\Qcon\dth{;}\Om_k')\\
        &\leq2\int_{\sp_+}D_\infty(\answ{Q_\phi}(\omega))\ \answ{d\H^2(\omega)}=\int_\sp D_\infty(\answ{Q_\phi}(\omega))\ \answ{d\H^2(\omega)}.
    \end{align*}
\end{proof}

\subsection{The Case $\lambda\in [0,\infty)$}
\label{ss:upII}

In this case we are not able to get an explicit equation for minimizers on a ray, so instead we use a similar method as in \cite{abl} where we approximate the minimizer by a step function convoluted with a smooth convolution kernel. 
\begin{proposition}\label{prop:upper_lambdafinite}
    If $\Qm$  minimizes $\exe$ with $\Qm=\answ{Q_\phi}$ on $\sp$ and\label{prop:up:fin}
    \begin{equation*}
        \frac{\eta}{\xi}\to\lambda\in [0,\infty),\quad\text{as}\quad\xi\dth{,\eta}\to0,
    \end{equation*}
    then
    \begin{equation*}
        \limsup_{\xi\dth{,\eta}\to0}\eta \exe(\Qm)\leq\int_\sp D_\lambda(\answ{Q_\phi}(\omega))\ \answ{d\H^2(\omega)}.
    \end{equation*}
\end{proposition}
\begin{proof}
    The goal of this proof will be to construct an equivariant sequence of maps which, in the limit $\xi\dth{,\eta}\to0$, attain the lower bound proven in Section 2, so we consider the map
    \begin{equation*}
        \overline{\Qc}(r,\theta,\phi)=R_\theta^T\Qc(r,\phi) R_\theta
        \, ,
    \end{equation*}
    where we will define $\Qc$ on the region $\Om_0:=\{(r,\phi):1<r<\infty,\ 0<\phi<\pi\}$. We introduce two parameters $h,\e>0$ to help with this construction, but their purposes will be discussed further as needed.
    
    \answ{Although we do this construction for boundary data $Q_\phi$, as noted in Remark \ref{rmk}, the same construction could be done for any $Q_b$ which is equivariant, Lipschitz, and which agrees locally with $Q_\phi$ near the poles. As a result, we will use the notation $Q_b$ at any points where the analysis would be identical for more general boundary data $Q_b$.}
    
    The strategy will again be to subdivide $\Om_0$ into smaller regions and define $\Qc$ piecewise on these smaller regions (see Figure \ref{fig:Ome0}). 
    \answ{In each region, we will construct a competitor and show that the only contribution to the energy in the limit $\xi,\eta\to 0$ comes from the region $\Omega_1$. More precisely} we will show that,
    \begin{equation*}
        \lim_{h\to 0}\left(\lim_{\e\to0}\left(\limsup_{\xi\dth{,\eta}\to0}\eta \exe(\overline{\Qc})\right)\right)\leq \int_\sp D_\lambda(\answ{Q_\phi}(\omega))\ \answ{d\H^2(\omega)}
        \, ,
    \end{equation*}
    which will imply that,
    \begin{equation*}
        \limsup_{\xi\dth{,\eta}\to0}\eta \exe(\Qm)\leq\int_\sp D_\lambda(\answ{Q_\phi}(\omega))\ \answ{d\H^2(\omega)}
        \, .
    \end{equation*}

    \begin{figure}[h]
    \begin{center}
    \begin{tikzpicture}[scale=1.5]

    \draw (0,1) arc (90:-90:1);
    \draw (0,1) -- (0,2.5);
    \draw (0,-1) -- (0,-2.5);
    \draw (45:1)-- ++(45:1);
    \draw (3*45/2:1)-- ++(3*45/2:0.5);
    \draw (90:1.5) arc (90:3*45/2:1.5);
    \draw (-3*45/2:1)-- ++(-3*45/2:0.5);
    \draw (-90:1.5) arc (-90:-3*45/2:1.5);
    \draw (-45:1)-- ++(-45:1);
    \draw (90:2) arc (90:-90:2);
    \draw[<->] (0,0.9) arc (90:45:.9);
    \draw (70:0.8) node [below] {$2h\eta$};
    \draw[<->] (40:1.05) -- ++(45:0.9);
    \draw (41:1.4) node [below right] {$2h\eta$};
    \draw (3,0) node {$\Omega_1$};
    \draw (1.5,0) node {$\Omega_2$};
    \draw (0.25,1.2) node {$\Omega_3^+$};
    \draw (0.9,1.4) node {$\Omega_4^+$};
    \draw (0.25,-1.2) node {$\Omega_3^-$};
    \draw (0.9,-1.4) node {$\Omega_4^-$};
\end{tikzpicture}
\end{center}
\caption{Subdivision of $\Om_0$ for the proof of Proposition \ref{prop:up:fin}}
\label{fig:Ome0}
\end{figure}
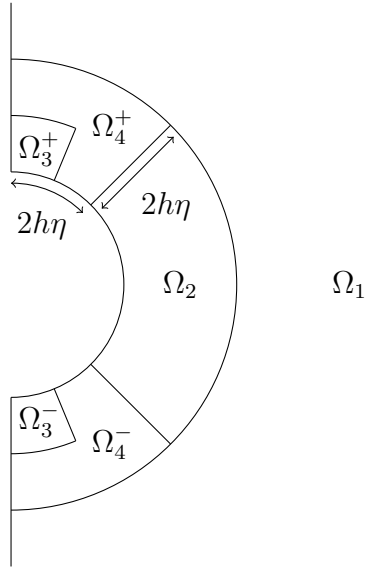
    \noindent \underline{\textit{Energy in $\Omega_1$}}: On this region, we will define the sequence of maps to be a convolution of a smooth kernel, $\varphi_\e$ of the angular coordinate $\phi$, and a function $Q^h$, which is a step function in $\phi$ and smooth in $r$. We want to consider the region,
    \begin{equation*}
        \Om_1:=\{(r,\phi):1+2h\eta<r<\infty,\ 0<\phi<\pi\}
    \end{equation*}
    and we begin by constructing a function $Q^h$ independent of $\xi$. While $H^1$ in $r$, the function $Q^h$ will be chosen to be piecewise constant w.r.t.\ $\phi$. To this goal, choose a partition
    \begin{equation*}
        0=\phi_0^h<\phi_1^h<\dots<\phi_{I_h}^h=\pi,\quad\text{with}\quad\phi_{i+1}^h-\phi_i^h<h
    \end{equation*}
    for all $i=1,\dots,I_h-1$. Then at each $i=1,\dots, I_h-2$, we can fix $Q_i^h\in Q_\infty+C_c^\infty((1,\infty);\sym)$ close enough in $H^1$ to the minimizer of $F_\lambda$ with boundary condition $Q_i^h(1)=\answ{Q_b}(0,\phi_i^h)$, so that
    \begin{equation}
        F_\lambda(Q_i^h)\leq D_\lambda(\answ{Q_b}(0,\phi_i^h))+h.\label{Dh}
    \end{equation}
    \answ{We then} define
    \begin{align*}
        Q^h(r,\phi)=\begin{cases}
            Q_i^h(r),&\phi\in[\phi_i^h,\phi_{i+1}^h),\quad i=1,\dots,I_h-2\\
            Q_\infty,&\phi\in(0,\phi_1^h)\cup[\phi_{I_h-1}^h,\pi).
        \end{cases}
    \end{align*}
    Take $\varphi_\e(\phi)=\e^{-1}\varphi(\phi/\e)$ to be a smooth kernel with $\supp\varphi=[-1,1]$, $\varphi\geq0$ and $\Vert\varphi\Vert_{L^1}=1$. 
    \answ{Finally, we} define $\Qc$ and $Q^{h,\e}$ as follows
    \begin{equation*}
        \Qc(r,\phi):=Q^{h,\e}(\tilde r,\phi):=\left(\varphi_\e*Q^h(\tilde{r},\cdot)\right)(\phi).
    \end{equation*}
    In this definition, we have
    \begin{equation}
        \tilde r=\frac{r-1-2h\eta}{\eta}+1\label{tilder}
        \, .
    \end{equation}
    Later on we will need to consider $r^2$ in terms of $r,h$ and $\eta$, so in order to simplify future calculations, we let
    \begin{align*}
        P(\tilde r,\eta, h):&=\frac{1}{\eta}(r^2-1)=\frac{1}{
        \eta}\big((\eta(\tilde r-1)+2h\eta+1)^2-1\big)
        \\
        &=\eta(\tilde r-1+2h)^2+2(\tilde r-1+2h).
    \end{align*}
    We note that
    \begin{equation}
        r^2=1+\eta P(\tilde r,\eta,h).\label{poly}
    \end{equation}
    As well, since $P$ is a polynomial in these three variables, if each of them is bounded, then $P$ is bounded.
    
   \answ{Next} we observe that $Q_i^h-Q_\infty$ is compactly supported for each $i$, so there exists a number $R(h)$ such that $\supp(Q_i^h-Q_\infty)\subset[1,R(h)]$ for all $i$. 
    Thus $Q^{h,\e}-Q_\infty$ is smooth and compactly supported on $[1,R(h)]\times(0,\pi)$. 
    We also have that $f$ and $g$ are both continuous, so there exists a constant $C(h,\e)>0$ such that
    \begin{equation*}
        \abs{\frac{\partial Q^{h,\e}}{\partial \tilde r}}^2,f(Q^{h,\e}),g(Q^{h,\e})\leq C(h,\e).
    \end{equation*}
    \answ{If} we take $\e<h/4$, which is justified since we will take $\e\to0$ before we take $h\to0$, then using the definition of $\varphi_\e$ and $Q^{h,\e}$ we get that, 
    \begin{equation*}
        Q^{h,\e}=Q_\infty\quad \text{on}\quad (0,\e]\cup[\pi-\e,\pi)
    \end{equation*}
    The final observations we make before addressing the energy are on the derivatives in the $\phi$ and $\theta$ directions. First note that from \cite{abl}, we have the bound
    \begin{equation}
        \abs{\frac{\partial Q}{\partial \theta}}^2\leq C|Q-Q_\infty|^2\label{theta}
    \end{equation}
    for any $Q\in Q_\infty+H^1(\Omega;\sym)$. Also, since $Q^{h,\e}-Q_\infty$ is compactly supported and smooth,
    \begin{equation*}
        \abs{Q^{h,\e}-Q_\infty}^2,\abs{\frac{\partial Q^{h,\e}}{\partial \phi}}^2<C(h,\e)
    \end{equation*}
    for some constant $C(h,\e)>0$. 
    We \answ{can then} consider the energy in this region:
    \begin{align}
        \eta \exe(\overline{\Qc}\dth{;}\Om_1')&=2\pi\eta\int_0^\pi\int_{1+2h\eta}^\infty \frac{r^2\sin\phi}{2}\abs{\frac{\partial \Qc}{\partial r}}^2+\frac{\sin\phi}{2}\abs{\frac{\partial \Qc}{\partial \phi}}+\frac{1}{2\sin\phi}\abs{\frac{\partial\overline{\Qc}}{\partial \theta}}\notag\\
        &\quad\quad+\frac{r^2\sin\phi}{\xi^2}f(\Qc)+\frac{r^2\sin\phi}{\eta^2}g(\Qc)\ drd\phi \label{separ} \\
        &=:A_r+A_\phi+A_\theta+A_f+A_g \nonumber
        \, .
    \end{align}
    To make this quantity easier to handle, we rename each term in a suggestive way by $A_r$, $A_\phi$, $A_\theta$, $A_f$, and $A_g$. Starting with $A_r$, we do a change of variables using \eqref{tilder} and the identity from \eqref{poly},
    \begin{align*}
        A_r&=2\pi\eta\int_0^\pi\int_{1+2h\eta}^\infty\frac{r^2\sin\phi}{2}\abs{\frac{\partial \Qc}{\partial r}}^2\ drd\phi\\
        &=2\pi\int_0^\pi\int_1^\infty\frac{(1+\eta P(\tilde r,\eta,h))\sin\phi}{2}\abs{\frac{\partial Q^{h,\e}}{\partial \tilde r}}^2\ d\tilde r d\phi.\label{rint}
    \end{align*}
    \answ{Next,} we separate the integral into two parts, and apply the bound from above,
    \begin{align*}
        A_r\leq2\pi\int_0^\pi\int_1^\infty\frac{\sin\phi}{2}\abs{\frac{\partial Q^{h,\e}}{\partial \tilde r}}^2\ d\tilde r d\phi+\pi\eta\int_0^\pi\int_1^{R(h)}P(\tilde r,\eta,h)C(h,\e)\ d\tilde rd\phi.
    \end{align*}
    For $\tilde r\leq R(h)$ and $h,\eta$ sufficiently small, $P$ is bounded by a constant, $C>0$, so the entire second term is bounded by $\eta C(h,\e)$, for some constant $C(h,\e)$. Thus it will disappear as $\xi\dth{,\eta}\to0$, for fixed $h,\e$. Also the first term, does not have any dependence on $\xi$, thus
    \begin{equation*}
        \limsup_{\xi\dth{,\eta}\to0}A_r\leq2\pi\int_0^\pi\int_1^\infty\frac{\sin\phi}{2}\abs{\frac{\partial Q^{h,\e}}{\partial \tilde r}}^2\ d\tilde r d\phi.
    \end{equation*}
    Next we consider $A_{\phi}$ and $A_\theta$ together, and using \eqref{theta} as well as the change of variables \eqref{tilder}, we have,
    \begin{align*}
        A_\phi+A_\theta\, &=\, 2\pi\eta\int_0^\pi\int_{1+2h\eta}^\infty\frac{\sin\phi}{2}\abs{\frac{\partial \Qc}{\partial \phi}}^2+\frac{1}{2\sin\phi}\abs{\frac{\partial\overline{\Qc}}{\partial \theta}}^2\ drd\phi\\
        &\leq\,  2\pi\eta^2\int_0^\pi\int_1^\infty\frac{1}{2}\abs{\frac{\partial Q^{h,\e}}{\partial \phi}}^2+\frac{C|Q^{h,\e}-Q_\infty|^2}{2\sin\phi} \ d\tilde rd\phi
    \end{align*}
    \answ{From} the above bounds, we can see that
    \begin{equation*}
        A_\phi+A_\theta\, \leq \, 2\pi\eta^2\int_\e^{\pi-\e}\int_1^{R(h)}C(h,\e)+\frac{C(h,\e)}{2\sin\e}\ d\tilde rd \phi.
    \end{equation*}
    Clearly each term is bounded and the integral is over a set of finite measure, so $A_\phi+A_\theta\to 0$ as $\xi\dth{,\eta}\to0$.
    For $A_f$, we want to isolate $\lambda^2 f(Q^{h,\e})$, so we do a change of variables \eqref{tilder} and introduce some new terms as follows,
    \begin{align*}
        A_f
        \ &= \ 2\pi\eta\int_0^\pi\int_{1+2h\eta}^\infty \frac{r^2\sin\phi}{\xi^2}f(\Qc)\ drd\phi\\
        &= 2\pi\int_0^\pi\int_1^\infty (1+\eta P(\tilde r,\eta,h))\frac{\eta^2}{\xi^2}f(Q^{h,\e})\sin\phi\ d\tilde rd\phi\\
        &=2\pi\int_0^\pi\int_1^\infty\left(\lambda^2+\left(\frac{\eta^2}{\xi^2}-\lambda^2\right)+\frac{\eta^3}{\xi^2}P(\tilde r,\eta,h)\right)f(Q^{h,\e})\sin\phi\ d\tilde rd\phi.
    \end{align*}
    We can then split this integral into two pieces, one that has no dependence on $\xi\dth{,\eta}$ and one which goes to zero as $\xi\dth{,\eta}\to0$. 
    \answ{Indeed,}
    \begin{align*}
        A_f\ \leq \ 2\pi\int_0^\pi\int_1^\infty \lambda^2f(Q^{h,\e})\sin\phi\ d\tilde rd\phi
        +C(h,\e)\int_0^\pi\int_1^{R(h)}\left(\frac{\eta^2}{\xi^2}-\lambda^2\right)+\frac{\eta^3}{\xi^2}\ d\tilde rd\phi
        \, ,
    \end{align*}
    where we used that $P$ is bounded since $\tilde r$ is bounded and both $\eta,h$ are small. As $\frac{\eta}{\xi}\to\lambda$ for $\xi\dth{,\eta}\to0$, the second integral disappears when $\xi\dth{,\eta}\to0$, and we are left with,
    \begin{equation*}
        \limsup_{\xi\dth{,\eta}\to0} A_f\leq 2\pi\int_0^\pi\int_1^\infty \lambda^2f(Q^{h,\e})\sin\phi\ d\tilde rd\phi.
    \end{equation*}
    Finally, we consider $A_g$,
    \begin{align*}
        A_g&=2\pi\eta\int_0^\pi\int_{1+2h\eta}^\infty\frac{r^2\sin\phi}{\eta^2}g(\Qc)\ drd\phi\\
        &=2\pi\int_0^\pi\int_1^\infty(1+\eta P(\tilde r,\eta,h))g(Q^{h,\e})\sin\phi\ d\tilde rd\phi.
    \end{align*}
    We again did the change of variable \eqref{tilder} and we use the same trick of splitting the integral into a piece that disappears and a piece that has no dependence on $\xi$.
    \begin{equation*}
        A_g
        \leq
        2\pi\int_0^\pi\int_1^\infty g(Q^{h,\e})\sin\phi\ d\tilde rd\phi
        +
        2\pi\eta\int_0^\pi\int_1^{R(h)}P(\tilde r,\eta,h)C(h,\e)\ d\tilde rd\phi
        \, ,
    \end{equation*}
    therefore,
    \begin{equation*}
        \limsup_{\xi\dth{,\eta}\to0}A_g\leq2\pi\int_0^\pi\int_1^\infty g(Q^{h,\e})\sin\phi\ d\tilde rd\phi.
    \end{equation*}
    We are now ready to put all the terms together, using \eqref{separ} to get the following bound on the energy,
    \begin{align*}
        \limsup_{\xi\dth{,\eta}\to0}\eta \exe(\overline{\Qc}\dth{;}\Om_1')&\leq2\pi\int_0^\pi\int_1^\infty \frac{1}{2}\abs{\frac{\partial Q^{h,\e}}{\partial\tilde r}}^2+\lambda^2f(Q^{h,\e})+g(Q^{h,\e})\ d\tilde r\ \sin\phi\ d\phi\\
        &=2\pi\int_0^\pi F_\lambda(Q^{h,\e}(\cdot,\phi))\sin\phi\ d\phi.
    \end{align*}
    Next we want to take $\e\to 0$, and we note that $Q^{h,\e}\to Q^h$ pointwise a.e., so by continuity of $f$ and $g$, $f(Q^{h,\e})\to f(Q^h)$ and $g(Q^{h,\e})\to g(Q^h)$ a.e. as $\e\to0$. As well,
    \begin{equation*}
        \frac{\partial Q^{h,\e}}{\partial r}=\varphi_\e*\frac{\partial Q^h}{\partial r}\to\frac{\partial Q^h}{\partial r}\quad\text{pointwise as}\quad\blue{\e}\to0,
    \end{equation*}
    thus,
    \begin{align*}
        \lim_{\e\to0}\left(\limsup_{\xi\dth{,\eta}\to0}\eta \exe(\overline{\Qc}\dth{;}\Om_1')\right)&\leq\lim_{\e\to0}2\pi\int_0^\pi F_\lambda(Q^{h,\e}(\cdot,\phi))\sin\phi\ d\phi\\
        &=2\pi\int_0^\pi F_\lambda(Q^h(\cdot,\phi))\sin\phi\ d\phi.
    \end{align*}
    Next we want to use that $Q\mapsto D_\lambda(Q)$ is locally Lipschitz for $Q\in\sym$, thus it is uniformly Lipschitz on the compact set of uniaxial $Q$-tensors. We note that the map $\phi\mapsto {Q_b}(0,\phi)$ is also Lipschitz, so their composition is Lipschitz with some constant $\L>0$. Using this, the piecewise definition of $Q^h$ and \eqref{Dh} we get,
    \begin{align*}
        \int_0^\pi F_\lambda(Q^h(\cdot,\phi))\sin\phi\ d\phi&=\sum_{i=1}^{I_h-2}\int_{\phi_i^h}^{\phi_{i+1}^h}F_\lambda(Q_i^h)\sin\phi\ d\phi\\
        &\leq \sum_{i=1}^{I_h-2}\int_{\phi_i^h}^{\phi_{i+1}^h} (D_\lambda(Q_b(0,\phi_i^h))+h)\sin\phi\ d\phi\\
        &\leq\sum_{i=1}^{I_h-2}\int_{\phi_i^h}^{\phi_{i+1}^h} (D_\lambda(Q_b(0,\phi))+(\L+1)h)\sin\phi\ d\phi
    \end{align*}
    \answ{Since} the integrand does not depend on $i$, we can combine the integrals and remove the $h$ term to get,
    \begin{equation*}
        \int_0^\pi F_\lambda(Q^h(\cdot,\phi))\sin\phi\ d\phi\leq\int_0^\pi D_\lambda({Q_b}(0,\phi)\sin\phi\ d\phi+(\L+1)\pi h.
    \end{equation*}
    Finally, we take $h\to 0$ and we see that,
    \begin{align*}
        \lim_{h\to0}\left(\lim_{\e\to0}\left(\limsup_{\xi\dth{,\eta}\to0}\eta \exe(\overline{\Qc}\dth{;}\Om_1')\right)\right)&\leq2\pi\int_0^\pi D_\lambda({Q_b}(0,\phi)\sin\phi\ d\phi\notag\\
        &=\int_\sp D_\lambda({Q_b}(\omega))\ \answ{d\H^2(\omega)}.
    \end{align*}
    Therefore it remains to show that the upper bound on the energy for all the other regions tends to zero.
    
    \noindent \underline{\textit{Energy in $\Omega_2$}}: On this region, we define our map $\Qc$ by a radial interpolation between $Q_b$ on the surface of the sphere and $Q^{h,\e}(1,\phi)$, as defined on $\Omega_1$, when $r=1+2h\eta$. We let,
    \[\Om_2:=\{(r,\phi):1<r<1+2h\eta,\ 2h\eta<\phi<\pi-2h\eta\},\]
    and we define $\Qc$ as follows,
    \begin{equation*}
        \Qc(r,\phi)=\left(\frac{r-1}{2h\eta}\right)Q^{h,\e}(1,\phi)+\left(1-\frac{r-1}{2h\eta}\right)Q_b(0,\phi).
    \end{equation*}
    It is clear that $Q^{h,\e}$ is uniformly bounded as we take $\xi\dth{,\eta}\to0$, since
    \begin{align*}
    |\Qc|\leq|Q^{h,\e}(1,\phi)|+|Q_b(0,\phi)|\leq|Q_b|_{L^\infty}|\varphi_\e|_{L^1}+|Q_b|_{L^\infty}=2\sqrt{\frac{2}{3}}.
    \end{align*}
    Using that $f$ and $g$ are continuous, it also holds that $|\Qc-Q_\infty|^2,f(\Qc),g(\Qc)\leq C$.
    Next we consider the derivative terms, beginning with the radial derivative,
    \begin{equation}
        \abs{\frac{\partial \Qc}{\partial r}}^2=\frac{1}{4h^2\eta^2}\abs{Q^{h,\e}(1,\phi)-Q_b(0,\phi)}^2.\label{r:der}
    \end{equation}
    If we consider $\phi\in[\phi_i^h,\phi_{i+1}^h)$, for $i=2,\dots,I_h-3$, then $Q^h(1,\phi\pm\e)=Q_b(0,\phi_i^h)$ or $Q^h(1,\phi\pm\e)=Q_b(0,\phi_{i\pm1}^h)$, when $\e$ is sufficiently small. In either case, for $-\e<\tau<\e$,
    \begin{equation*}
        |Q^h(1,\phi-\tau)-Q_b(0,\phi)|\leq 2h\L_b
        \, ,
    \end{equation*}
    where $\L_b$ is the Lipschitz constant of $\phi\mapsto Q_b(0,\phi)$.
    This will be useful in bounding the derivative as
    \begin{align*}
        \abs{\frac{\partial \Qc}{\partial r}}^2&\leq\frac{1}{4h^2\eta^2}\abs{\int_{-\e}^\e \varphi_\e(\tau)Q^h(1,\phi-\tau)\ d\tau-Q_b(0,\phi)\int_{-\e}^\e\varphi_\e(\tau)\ d\tau}^2\notag\\
        &\leq\frac{1}{4h^2\eta^2}\abs{\int_{-\e}^\e \varphi_\e(\tau)|Q^h(1,\phi-\tau)-Q_b(0,\phi)|\ d\tau}^2\leq\frac{\L_b^2}{\eta^2}
        \, ,
    \end{align*}
    which holds for all $\phi\in(2h\eta+2h,\pi-2h\eta-2h)$. If we consider $\phi$ outside of this range, then we simply use the bound,
    \begin{equation}
        \abs{Q^h(1-\phi,\tau)-Q_b(0,\phi)}\leq 2|Q_b|_{L^\infty}\leq C\label{bd:hconst}
    \end{equation}
    for some constant $C>0$. 
    \answ{For} the radial derivative term in the energy and for $\phi\in(2h\eta+2h,\pi-2h\eta-2h)$ it holds that
    \begin{align*}
        \int_{2h\eta+2h}^{\pi-2h\eta-2h}\int_1^{1+2h\eta}\frac{r^2\sin\phi}{2}\abs{\frac{\partial \Qc}{\partial r}}^2\ drd\phi
        \,\leq\,
        \int_0^\pi\int_1^{1+2h\eta}\frac{(1+2h\eta)^2\L_b^2}{2\eta^2}\ drd\phi
        \, \leq\, 
        \frac{C h}{\eta}
    \end{align*}
    for some constant $C>0$, assuming $\eta$ is sufficiently small. 
    For the remaining values of $\phi$, let $A=(2h\eta,2h\eta+2h)\cup(\pi-2h\eta-2h,\pi-2h\eta)$, then using~\eqref{bd:hconst}
    \begin{align*}
        \int_A\int_1^{1+2h\eta}\frac{r^2\sin\phi}{2}\abs{\frac{\partial \Qc}{\partial r}}^2\ drd\phi&\leq\int_A\int_1^{1+2h\eta}\frac{C(1+2h\eta)^2\sin(2h\eta+2h)}{h^2\eta^2}\ drd\phi\notag\\
        &\leq 4h(2h\eta)\left[\frac{C(1+2h\eta)^2(\eta+1)}{h\eta^2}\right]\leq\frac{Ch}{\eta}
    \end{align*}
    for some constant $C>0$. So the whole radial derivative term has the same bound. Next we consider the $\phi-$derivative.
    \begin{align*}
        \abs{\frac{\partial \Qc}{\partial \phi}}^2&=\abs{\left(\frac{r-1}{2h\eta}\right)\frac{Q^{h,\e}}{\partial \phi}(1,\phi)+\left(1-\frac{r-1}{2h\eta}\right)\frac{\partial Q_b}{\partial \phi}(0,\phi)}^2\\
        &\leq\abs{\frac{\partial Q^{h,\e}}{\partial \phi}(1,\cdot)}_{L^\infty}^2+2\abs{\frac{\partial Q^{h,\e}}{\partial \phi}(1,\cdot)}_{L^\infty}\abs{\frac{\partial Q_b}{\partial \phi}(0,\cdot)}_{L^\infty}+\abs{\frac{\partial Q_b}{\partial \phi}(0,\cdot)}_{L^\infty}^2\,\leq\, C(h,\e)\label{bd:phi}
        \, ,
    \end{align*}
    where $C(h,\e)>0$ is some constant which depends on $h$ and $\e$, since $Q_b$ is Lipschitz and $Q^{h,\e}(1,\cdot)$ is smooth. 
    We \answ{can then} bound the entire energy on $\Om_2'$. 
    Using the above bounds on each term, we see that
    \begin{multline*}
        \eta \exe(\overline{\Qc}\dth{;}\Om_2')\leq2\pi\eta\left(\frac{Ch}{\eta}+\int_{2h\eta}^{\pi-2h\eta}\int_1^{1+2h\eta}\frac{C(h,\e)\sin\phi}{2}+\frac{C}{\sin\phi}+\frac{Cr^2\sin\phi}{\xi^2}\right.\\
        \left.+\frac{Cr^2\sin\phi}{\eta^2}\ drd\phi\right).
    \end{multline*}
    We then use that $\sin\phi\geq\sin(2h\eta)\geq h\eta$ to get,
    \begin{align*}
        \eta \exe(\overline{\Qc}\dth{;}\Om_2')&\leq Ch+2\pi^2\eta(2h\eta)\left[C(h,\e)+\frac{C}{h\eta}+\frac{C}{\xi^2}+\frac{C}{\eta^2}\right]\\
        &\leq Ch+C(h,\e)h\eta^2+C\eta+\frac{Ch\eta^2}{\xi^2}.
    \end{align*}
    Next we take $\xi\dth{,\eta}\to0$ and use that $\frac{\eta}{\xi}\to\lambda$ to get
    \begin{equation*}
        \limsup_{\xi\dth{,\eta}\to0}\eta \exe(\overline{\Qc}\dth{;}\Om_2')\leq Ch+Ch\lambda^2.
    \end{equation*}
    Since there is no dependence on $\e$ anymore, and each term has a factor of $h$,
    \begin{equation*}
        \lim_{h\to0}\left(\lim_{\e\to0}\left(\limsup_{\xi\dth{,\eta}\to0}\eta \exe(\overline{\Qc}\dth{;}\Om_2')\right)\right)=0.
    \end{equation*}
    
    \noindent \underline{\textit{Energy in $\Omega_3$}}: For $\Om_3$, we will take $\Qc$ to be uniaxial and we will use almost the same definition as in $\Om_3$ for $\lambda=\infty$. \answ{In doing this, we are imposing that $Q_b=Q_\phi$ on the piece of $\partial \Om_3$ which coincides with $\partial B$. As a result, we rely on $Q_b=Q_\phi$ in this step of the construction.} Define $\Omega_3$ in two pieces, $\Om_3:=\Omega_3^+\cup \Om_3^-$, where
    \begin{align*}
        \Om_3^+&:=\{(r,\phi):1<r<1+h\eta,\ 0<\phi<h\eta\},\\
        \Om_3^-&:=\{(r,\phi):1<r<1+h\eta,\ \pi-h\eta>\phi<\pi\}
        \, ,
    \end{align*}
    then we define $\Qc$ on $\Om_3^+$ by,
    \begin{equation*}
        \Qc(r,\phi)
        \ = \
        R_\phi m\left(\frac{r-1}{h\eta},\frac{\phi}{h\eta}\right)\otimes R_\phi m\left(\frac{r-1}{h\eta},\frac{\phi}{h\eta}\right)-\frac{1}{3}I
        \, ,
    \end{equation*}
    where $R_\phi$ and $m$ are defined from \eqref{sing1} and \eqref{sing2}. We note the only difference is that rather than defining it on a region of size $\eta^2$, it is on a region of size $h^2\eta^2$, so by a simple change of variables $q=h\eta$, we can do all the previous computations and replace $\eta$ by $q$. 
    Then, from \eqref{sing:en}, we get the bound
    \begin{equation*}
        \eta \exe(\overline{\Qc}\dth{;}(\Om_3^+)')\leq 2\pi C\eta q
    \end{equation*}
    and since $q\to0$ as $\eta\to 0$, we see that $\eta \exe(\overline{\Qc}\dth{;}(\Om_3^+)')\to 0$ as $\xi\dth{,\eta}\to0$.
    Again we get the same bound on the energy in $\Om_3^-$ by a reflection about the $\{x_3=0\}$ plane, therefore,
    \begin{equation*}
        \limsup_{\xi\dth{,\eta}\to0}\eta \exe(\overline{\Qc}\dth{;}\Om_3')=0.
    \end{equation*}
    Since the right hand side has no dependence on $\e$ or $h$, we can take the limit $\e,h\to 0$.

    \noindent \underline{\textit{Energy in $\Omega_4$}}: We will approach $\Om_4$ in a very similar way to the $\lambda=\infty$ case, with a Lipschitz extension. Here we treat each component, $Q_{ij}$, of the $Q$-tensor as its own function and do a Lipschitz extension of each, but we do this carefully as to ensure that the extension is still symmetric and traceless. Let $\Om_4:=\Om_4^+\cup \Om_4^-$, where
    \begin{align*}
        \Om_4^+&:=\{(r,\phi):1<r<1+2h\eta,\ 0<\phi<2h\eta\}\setminus \overline{\Om_3^+},\\
        \Om_4^-&:=\{(r,\phi):1<r<1+2h\eta,\ \pi-2h\eta<\phi<\pi\}\setminus\overline{\Om_3^-}.
    \end{align*}
    We will begin by working on the upper bound for the energy on $\Om_4^+$, but it will be clear later that by the symmetry of the construction, we can get the same bound for the energy on $\Om_4^-$. First we will find a Lipschitz constant for $\Qc$ on the boundary by finding an upper bound on the derivatives along each piece of the boundary. We use the same naming of boundary components as in the $\lambda=\infty$ case, labelling the sides $\ell_1,\dots,\ell_6$. On $\ell_1$, $\Qc=Q_b$, so it has Lipschitz constant $\L_b$ \answ{which has no dependence on $\eta$}. Next we have $\ell_2$ where using \eqref{sing:r} we get,
    \begin{equation*}
        \abs{\frac{\partial \Qc}{\partial r}}^2=2\abs{\frac{\partial (R_\phi m)}{\partial r}}^2\leq \frac{C}{h^2\eta^2}\left(\frac{\tau^2}{|(s,\tau)|^4}\right)\leq\frac{C}{h^2\eta^2}
    \end{equation*}
    since we are on the line where $\tau=1$. Similarly on $\ell_3$ we use \eqref{sing:phi} to see that
    \begin{equation*}
        \abs{\frac{\partial \Qc}{\partial \phi}}^2=2\abs{\frac{\partial (R_\phi m)}{\partial \phi}}^2\leq \frac{C}{h^2\eta^2}\left(\frac{s^2}{|(s,\tau)|^4}\right)\leq\frac{C}{h^2\eta^2}
    \end{equation*}
    as we are evaluating at $s=1$. 
    Along $\ell_4$ and $\ell_5$, $\Qc=Q_\infty$ for sufficiently small $\eta$, so the Lipschitz constant is zero. Finally on $\ell_6$, using \eqref{r:der} and \eqref{bd:hconst} we have,
    \begin{equation*}
        \abs{\frac{\partial \Qc}{\partial r}}^2\leq\frac{C}{h^2\eta^2}|Q^{h,\e}(1,\phi)-Q_b(0,\phi)|\leq \frac{C}{h^2\eta^2}
    \end{equation*}
    since $Q^{h,\e}=Q_\infty$ here. Therefore $\Qc$ has the Lipschitz constant $\tilde{\L}=\frac{C}{h\eta}$ for some constant $C>0$, so each component function of the $Q$-tensor, $\Qc$ also has the same Lipschitz constant. We take each $Q_{ij}$ above the main diagonal and define a Lipschitz extension for them on $\Om_4^+$ and we do the same for $Q_{11}$ and $Q_{22}$, then we let $Q_{33}=-Q_{11}-Q_{22}$ to maintain the traceless condition and define $Q_{ji}=Q_{ij}$ so that the extension is still symmetric. 
    \answ{Then} $\Qc$ is a Lipschitz function on $\Om_4^+$ and the derivatives are bounded as follows
    \begin{equation*}
        \abs{\frac{\partial \Qc}{\partial r}}^2,\abs{\frac{\partial \Qc}{\partial \phi}}^2\leq \frac{C}{h^2\eta^2}.
    \end{equation*}
    Furthermore, to bound the $\theta-$derivative term, we can see that, \begin{equation*}
        |\Qc-Q_\infty|^2=|\Qc(r,\phi)-\Qc(r,0)|^2\leq \tilde{\L}^2\phi^2.
    \end{equation*}
    Finally $\Qc$ is bounded independently of $h$ and $\eta$ on the boundary of $\Om_4^+$, and since the diameter of $\Om_4^+$ is proportional to $h\eta$ while the Lipschitz constant is proportional to $h^{-1}\eta^{-1}$, the function $\Qc$ is bounded by a constant on $\Om_4^+$. Therefore $f(\Qc),g(\Qc)\leq C$ for some $C>0$. Now we can compute the energy on the region to be
    \begin{align*}
        \eta \exe(\overline{\Qc}\dth{;}(\Om_4^+)')&=2\pi\eta\iint_{\Om_4^+}\frac{r^2\sin\phi}{2}\abs{\frac{\partial \Qc}{\partial r}}^2+\frac{\sin\phi}{2}\abs{\frac{\partial \Qc}{\partial \phi}}^2+\frac{C|\Qc-Q_\infty|^2}{\sin\phi}\\
        &\qquad\qquad\qquad\qquad+\frac{r^2\sin\phi}{\xi^2}f(\Qc)+\frac{r^2\sin\phi}{\eta^2}g(\Qc)\ drd\phi\\
        &\leq 2\pi\eta\iint_{\Om_4^+}\frac{C}{h^2\eta^2}+\frac{C\phi^2}{\sin\phi}+\frac{C}{\xi^2}+\frac{C}{\eta^2}\ drd\phi\\
        &\leq C\eta(h^2\eta^2)\left(\frac{C}{h^2\eta^2}+Ch\eta+\frac{C}{\xi^2}+\frac{C}{\eta^2}\right)\\
        &=C\eta+Ch^3\eta^4+\frac{Ch^2\eta^3}{\xi^2}+Ch^2\eta
        \, .
    \end{align*}
    We notice that the same bounds all hold on $\Om_4^-$.
    Thus, when taking $\xi\dth{,\eta}\to0$, we see that
    \begin{align*}
        \limsup_{\xi\dth{,\eta}\to0}\eta \exe(\overline{\Qc}\dth{;}\Om_4')=0.
    \end{align*}
    Combining the bounds of all regions $\Omega_1,\ldots,\Omega_4$, it follows that 
    \begin{equation*}
        \limsup_{\xi\dth{,\eta}\to0}\eta \exe(\Qm)\leq\lim_{h\to0}\left(\lim_{\e\to0}\left(\limsup_{\xi\dth{,\eta}\to0}\eta \exe(\overline{\Qc})\right)\right)\leq\int_{\sp}D_\lambda(Q_b(\omega))\ \answ{d\H^2(\omega)}.
    \end{equation*}
    \answ{Recalling that $Q_b=Q_\phi$ concludes the proof of Proposition~\ref{prop:upper_lambdafinite}.}
\end{proof}

%

\paragraph{Disclosure statement.}
The authors report there are no competing interests to declare.

\paragraph{Data availability statement.} No data set associated with the paper.

\bibliographystyle{acm}

\end{document}